\documentclass[11pt]{amsart}

\usepackage{amsmath}
\usepackage{amsfonts}
\usepackage{amssymb} 
\usepackage{amsthm} 
\usepackage{mathptmx}
\usepackage{mathrsfs}
\usepackage{latexsym}
\usepackage{times}
\usepackage[isolatin]{inputenc}

\DeclareMathAlphabet{\mathpzc}{OT1}{pzc}{m}{it}

\newtheorem{thm}{Theorem}[subsection]
\newtheorem{lem}[thm]{Lemma}
\newtheorem{prop}[thm]{Proposition}
\newtheorem{cor}[thm]{Corollary}

\theoremstyle{definition}
\newtheorem{defn}[thm]{Definition}
\newtheorem{ex}[thm]{Example}

\theoremstyle{remark}
\newtheorem{rem}[thm]{Remark}

\newcommand{\gl}{\mathfrak{gl}}
\newcommand{\sk}{\mathfrak{s}}
\newcommand{\slk}{\mathfrak{sl}}

\newcommand{\zk}{\mathfrak{z}}
\newcommand{\Sk}{\mathfrak{S}}
\newcommand{\g}{\mathfrak{g}}

\newcommand{\hk}{\mathfrak{h}}
\newcommand{\lk}{\mathfrak{l}}
\newcommand{\ik}{\mathfrak{i}}

\newcommand\Fc{\mathcal F}

\newcommand\Dc{\mathcal D}
\newcommand\Ac{\mathcal A}
\newcommand\Ec{\mathcal E}
\newcommand\Hc{\mathcal H}
\newcommand\Pc{\mathcal P}
\newcommand\Cc{\mathcal C}
\newcommand\Mcc{\mathcal M}
\newcommand\Rc{\mathcal R}
\newcommand\Sc{\mathcal S}
\newcommand\Wc{\mathcal W}

\newcommand{\Dt}{\mathtt{\mathbf D}}
\newcommand{\Ft}{\mathtt{\mathbf F}}
\newcommand{\Hht}{\mathtt{H}}

\newcommand\RR{\mathbb R}
\newcommand\CC{\mathbb C}
\newcommand\NN{\mathbb N}
\newcommand\ZZ{\mathbb Z}

\newcommand{\Tr}{\operatorname{Tr}}

\newcommand\gd{\g^*}

\newcommand{\Cliff}{\operatorname{Cliff}}
\newcommand{\So}{\operatorname{S}}
\newcommand{\Ao}{\operatorname{A}}

\newcommand{\Zo}{\operatorname{Z}}
\newcommand{\io}{\operatorname{\iota}}

\newcommand{\ad}{\operatorname{ad}}
\newcommand{\im}{\operatorname{Im}}

\renewcommand\dfrac{\displaystyle \frac}
\newcommand\ds\displaystyle

\newcommand\e{\varepsilon}

\newcommand{\End}{\operatorname{End}}
\newcommand{\Der}{\operatorname{Der}}
\newcommand{\ext}{\operatorname{Ext}}
\newcommand{\sym}{\operatorname{Sym}}
\newcommand{\Id}{\operatorname{Id}}
\renewcommand\hat\widehat 
\renewcommand\tilde\widetilde
\newcommand\Ht{\tilde{H}} 
\newcommand{\spa}{\operatorname{span}}

\newcommand{\gla}{$\mathbf{gla}$ }
\newcommand{\GPB}{\operatorname{GPB}}

\newcommand{\rank}{\operatorname{rank}}

\newcommand\Wedge{\bigwedge} 
\newcommand\extg{\Wedge \g}
\newcommand\extgp{\Wedge_+ \g}
\newcommand\extgq{\Wedge_Q \g}
\newcommand\extgg{(\Wedge \g)^\g}
\newcommand\adpo{\ad_{\tt P}} 

\newcommand{\dss}{\displaystyle}

\begin{document}
\title[New applications of graded Lie algebras]{New applications of
  graded Lie algebras to Lie algebras, generalized Lie algebras and
  cohomology}

\author{Georges Pinczon, Rosane Ushirobira}

\address{Institut de Math\'ematiques de Bourgogne, Universit\'e de
  Bourgogne, B.P. 47870, F-21078 Dijon Cedex, France}

\email{gpinczon, rosane@u-bourgogne.fr}



\date{\today}

\begin{abstract} 
We give new applications of graded Lie algebras to: identities of
standard polynomials, deformation theory of quadratic Lie algebras,
cyclic cohomology of quadratic Lie algebras, $2k$-Lie algebras,
generalized Poisson brackets and so on.
\end{abstract}

\maketitle

\section{Introduction}

Graded Lie algebras ({\bf gla}) are commonly used in many areas of
Mathematics and Physics. One of the reasons is that they offer a very
convenient framework for the development of theories such as
Cohomology Theory, Deformation Theory, among others, very often
avoiding heavy computations. But the fundamental reason, developed
with the work of M. Gerstenhaber and others, is that the \gla notion
allows to endow with a structure, objects that a priori had none,
providing a new and efficient material to study these objects. The aim
of this paper is to give some new applications of classical well-known
\gla related to Deformation Theory.

Let us start with some notations: $\g$ will be a complex vector space,
$\extg$ the Grassmann algebra of $\g$, that is, the algebra of skew
multilinear forms on $\g$, with the wedge product. When $\g$ is
finite-dimensional, one has $\extg = \ext(\gd)$, where $\ext(\gd)$
denotes the exterior algebra of the dual space $\gd$. However, when
$\g$ is not finite dimensional, the strict inclusion $\ext(\gd)
\subset \extg$ holds. A {\em quadratic} vector space is a vector space
endowed with a non degenerate symmetric bilinear form. In the case of
a {\em quadratic Lie algebra} this bilinear form has to be
invariant. A theory of finite dimensional quadratic Lie algebras,
based on the notion of double extension, was developed in \cite{Santa,
  Medina} following Kac's arguments \cite{Kac}. In this paper, we
shall present another interpretation based on the concept of super
Poisson bracket.

The \gla we shall use here are:

\begin{enumerate}

\item Gerstenhaber's graded Lie algebras $\Mcc(\g)$, related to
  associative algebra structures on $\g$ (see Section \ref{Section1}).

\item Gerstenhaber-Nijenhuis's graded Lie algebras $\Mcc_a (\g)$,
  related to Lie algebra structures on $\g$ (see Section
  \ref{Section1}).

\item The graded Lie algebra $\Dc (\g)$ of derivations of the
  Grassmann algebra $\extg$ (often called $W(n)$ when $n= \dim \g$,
  see Section \ref{Section2}).

\item Assuming that $\g$ is finite dimensional, the graded Lie algebra
  $\Wc (\g)$ of skew symmetric polynomial multivectors on $\gd$ with
  the Schouten bracket (see Section \ref{Sections}).

\item Given a quadratic finite dimensional space $\g$, the super
  Poisson graded Lie algebra structure on the Grassmann algebra
  $\extg$ (see Section \ref{Section4}) and the superalgebra $\Hc(\g)$
  of Hamiltonian derivations of $\extg$.

\end{enumerate}

For (1) and (2), we refer to \cite{Gerst}, for (3) to \cite{Scheu},
for (4) to Koszul's presentation \cite{Koszul2} (though \cite{Gerst}
could be convenient as well). For (5), though it is a known algebra,
we have no references, probably because of the lack of applications up
to now (we shall show, e.g. in Sections 5 to 8, that there are some
natural and interesting ones). Since we want to fix our conventions
and notations and since we do not wish to address the present work to
\gla experts only, we give an introduction to all the above \gla,
recalling the main properties that will be used all along this paper.

Section \ref{Section1} is a review of $\Mcc (\g)$ and $\Mcc_a(\g)$. We
conclude the Section with a notion of Generalized Lie Algebras
structures, that we call $2k$-{\em Lie algebras}, namely the elements
$F$ in $\Mcc_a^{2k}(\g)$ that satisfy $[F,F]_a =0$. Such structures
are introduced in \cite{Filip} and many other papers
(e.g. \cite{Per}), under various names.

In Section \ref{Section2}, we recall how to go from $\Mcc_a(\g)$ to
$\Dc(\g)$, an operation that can be translated as going from a
structure to its cohomology, as we shall now explain. The argument is
given by (\ref{2.1.4}): there exists a one to one \gla homomorphism
from $\Mcc_a(\g)[1]$ to $\Dc(\g)$, which turns out to be an
isomorphism when $\g$ is finite dimensional. So given a $2k$-Lie
algebra structure on $\g$, there is an associated derivation $D$ of
$\extg$, and the (generalized) Jacobi identity $[F,F]_a =0$ is
equivalent to $D^2 =0$, so that $D$ defines a cohomology complex
(\ref{2.3}). This is well known for Lie algebras since the
corresponding complex is the Chevalley complex of trivial
cohomology. The existence of a cohomology complex for a $2k$-Lie
algebra was pointed out (without the \gla interpretation), e.g. in
\cite{Per}. We then recall the definition and properties of the
Schouten bracket for a finite dimensional $\g$. As in \cite{Per}, we
define a {\em Generalized Poisson Bracket} ($\GPB$) as an element $W$
of $\Wc^{2k} (\g)$ satisfying $[W,W]_S = 0$ (\ref{3.1.1}), the obvious
generalization of the classical definition of a Poisson bracket. We
show that there exists a one to one \gla homomorphism from $\Dc(\g)$
into $\Wc(\g)[1]$ (\ref{3.2.2}), so that any $2k$-Lie algebra
structure on $\g$ has an associated $\GPB$, generalizing the classical
Lie-Kostant-Kirillov bracket associated to a Lie algebra.

We apply the results of Sections \ref{Section1} and \ref{Section2} to
standard polynomials in Section \ref{Sections}. Standard polynomials
$\Ac_k$ ($k \geq 0$) on an associative algebra $\g$, appear in the PI
algebras theory (see \cite{Jacobson}) and also in cohomology theory
(for instance, the cohomology of $\gl(n)$ is $\ext[a_1, a_3, \dots,
  a_{2n-1}]$ where $a_k = \Tr(\Ac_k)$, and the cohomology of
$\gl(\infty)$ is $\ext[a_1, a_3, \dots]$ \cite{Fuchs}). We show that
there exist two different structures on the space $\Ac= \spa \{ \Ac_k
\mid k \geq 0 \}$, both with interesting consequences. The first one
comes from the Gerstenhaber bracket of $\Mcc_a(\g)$: we compute
explicitly $[\Ac_k, \Ac_{k'}]_a$, and it results that $\Ac$ is a
subalgebra of the \gla $\Mcc_a(\g)$ (\ref{s.2.2}). Since $[\Ac_{2k},
  \Ac_{2k}]_a = 0$, any even standard polynomial define a $2k$-Lie
algebra structure on $\g$ (\ref{s.2.5}). Moreover, $\Ac_{2k}$ is a
coboundary (an invariant one) of the adjoint cohomology of the Lie
algebra $\Ac_2$ defined by the associative algebra $\g$.  The second
product, denoted by $\times$, is associative and is defined on
$\Mcc_a(\g)$ using both the wedge product on $\extg$ and the initial
product on $\g$ (a priori, non commutative). So, it is rather a
surprise to find that $\Ac$ is an Abelian algebra for $\times$, and in
fact a very simple one, since $\Ac_k = (\Ac_1)^{\times k}$, $\forall
k$ (\ref{s.3.4}). For instance, for $\g = \gl(n)$, $\Ac$ with its
$\times$-product is isomorphic to $\CC[x]/ x^{2n}$, since $\Ac_k = 0$,
$\forall k \geq 2n$ (the Amitsur-Levitzki theorem \cite{Amitsur,
  Kostant1}). From identities (\ref{s.3.4}), one deduces some
classical well-known identities of standard polynomials
(e.g. $\Ac_{2k} = (\Ac_2)^{\times k}$, $\forall k$, usually proved by
hand). When $\g$ has a trace, we prove that $\Tr([F,G]_\times) = 0$,
for all $F$, $G \in \Mcc_a(\g)$ (\ref{s.3.5}), and then (keeping the
notation $a_k = \Tr(\Ac_k)$), that $a_{2k} = 0$, $\forall k > 0$, and
that $a_{2k+1}$ is an invariant Lie algebra cocycle (\ref{s.3.6}). To
conclude Section \ref{Sections}, we compute the cohomology of the Lie
algebra $\g$ of finite rank operators in an infinite dimensional
space. Obviously, $\gl(\infty) \subset \g$, but this inclusion is
strict. Our result is $H^\star (\g) = \ext[a_1, a_3, \dots]$
(\ref{s.8}), so the above inclusion induces an isomorphism in
cohomology.

The first part of Section \ref{Section4} is devoted to the
construction of the super Poisson bracket defined on $\extg$, when
$\g$ is a finite dimensional quadratic vector space. We follow a
deformation argument as in \cite{Kostant2}: the Clifford algebra
$\Cliff (\g^*)$ can be seen as a quantization of the algebra $\extg$
of skew polynomials, similarly to the classical Moyal quantization of
polynomials by the Weyl algebra. In \ref{4.1}, we introduce some
formulas for the construction of the Clifford algebra that are
convenient since they easily provide a transparent explicit formula
for the deformed product (\ref{4.1.6}), with leading term the super
Poisson bracket, explicitly computed in (\ref{4.2.1}). The relation
with the superalgebra $\tilde{H}(n)$ \cite{Scheu} is given in
(\ref{4.1.3}), and a Moyal type formula is obtained (\ref{4.2.2}) (an
equivalent formula without the super Poisson bracket can be found in
\cite{Kostant2}). In the second part of Section \ref{Section4}, we use
the \gla $\extg$ and the super Poisson bracket to study quadratic Lie
algebras.
We obtain that quadratic Lie algebra structures on $\g$ with bilinear
form $B$ are in one to one correspondence with elements $I$ in
$\Wedge^3 \g$ satisfying $\{I, I \} =0$; more precisely, $I(X,Y,Z) =
B([X,Y], Z)$, $\forall X$, $Y$, $Z \in \g$, and the differential
$\partial$ of $\extg$ is $\partial = - \frac12 \ \adpo(I)$
(\ref{4.4.2}, \ref{4.5.3}). We prove that any quadratic deformation of
a quadratic Lie algebra is equivalent to a deformation with unchanged
invariant bilinear form (\ref{4.4.1}), and finally, we propose a \gla
framework well adapted to the deformation theory of quadratic Lie
algebras (\ref{4.6.2}).

We use the results of Section \ref{Section4} in Section
\ref{Section4bis} to give a complete description of finite dimensional
elementary quadratic Lie algebras, i.e. those with decomposable
associated element $I$ in $\Wedge^3 \g$ (\ref{4.1.3ele}). We first
give a simple characterization (\ref{4bis.2.2}): a non Abelian
quadratic Lie algebra $\g$ is elementary if and only if $\dim ([\g,
  \g]) = 3$. We then show that any non Abelian quadratic Lie algebra
reduces, up to a central factor, to a quadratic Lie algebra with
totally isotropic center (\ref{4bis.4.1}); the property of being
elementary is preserved under the reduction. This reduces the problem
of finding all elementary non Abelian quadratic Lie algebras to
algebras of dimension 3 to 6 (\ref{4bis.4.2}), that we completely
describe in (\ref{4bis.5}) and (\ref{4bis.4.3}). Some remarks: as we
show in (\ref{4bis.2.3}), if $\g$ is an elementary quadratic Lie
algebra, all coadjoint orbits have dimension at most 2. Now, a
classification of Lie algebras whose coadjoint orbits are of dimension
less than 2, is given in \cite{Arnal}, and the proof, using classical
Lie algebra theory is not at all trivial. With some effort, one could
probably find directly in the classification of \cite{Arnal}, which
algebras are quadratic and which are not. Our geometric flavored proof
is completely different, using essentially elementary properties of
quadratic forms.

In Section \ref{Section5}, we study cyclic cohomology of quadratic Lie
algebras.  Given a quadratic vector space $\g$, we use $\g$-valued
cochains (rather than $\gd$-valued, by analogy to the associative case
\cite{Connes}) to define cyclic cochains (\ref{5.1.3}) (both notions
are equivalent when $\g$ is finite dimensional). Thanks to this
definition, we can use the Gerstenhaber bracket of $\Mcc_a (\g)$ and
we show that cyclic cochains are well behaved with respect to this
bracket: the space $\Cc_c(\g)$ of cyclic cochains is a subalgebra of
the \gla $\Mcc_a(\g)$ (\ref{5.1.4}) and if $\g$ is a Lie algebra,
$\Cc_c (\g)$ is a subcomplex of the adjoint cohomology complex $\Mcc_a
(\g)$ (\ref{5.4.3}); we define the cyclic cohomology $H_c^* (\g)$ as
the cohomology of this subcomplex (\ref{5.4.4}). There is a natural
one to one map from $\Cc_c(\g)$ into $\extgq = \extg / \CC$
(\ref{5.1.4}) which induces a map from $H^*_c(\g)$ into $H_Q^*(\g) =
H(\g) / \CC$. When $\g$ is finite dimensional, $\extgq$ is a \gla for
the (quotient) Poisson bracket, isomorphic to $\Hc(\g)$, and there is
an induced \gla structure on $H_Q(\g)$. We show that there is a \gla
isomorphism from $\Cc_c(\g)$ onto $\extgq$ (\ref{5.3.2bis}), and from
$H^*_c(\g)$ onto $H^*_Q(\g)$ (\ref{5.4.6bis}).  We also introduce a
wedge product on $\Cc_c(\g)$, and on $H^*_c(\g)$ (\ref{5.2.1},
\ref{5.2.3}) which proves to be useful to describe $H_c^*(\g)$
(\ref{5.4.6bis}). When $\g$ is not finite dimensional, the isomorphism
between $H^*_c(\g)$ and $H^*_Q(\g)$ is no longer true: we give an
example where the natural map is neither one to one, nor onto
(\ref{6.4.10}, \ref{6.4.11}). So the cyclic cohomology $H^*_c(\g)$ can
have its own life, independently of the reduced cohomology
$H^*_Q(\g)$.

Section \ref{Section6} starts with the study of invariant cyclic
cochains in the case of a finite dimensional quadratic Lie algebra. We
first prove that any invariant cyclic cochain is a cocycle
(\ref{6.1.5}). When $\g$ is reductive, we demonstrate that the
inclusion of invariants cyclic cochains into cocycles induces an
isomorphism in cohomology (\ref{6.1.5}), so that $H_c^* (\g) \simeq
\Cc_c(\g)^\g$. Assuming that $\g$ is a semisimple Lie algebra, we
prove:

{\centerline{{\em If $I, I' \in \extgg$, then $\{ I, I' \} = 0$.}
  (\ref{6.2.1})}}

As a corollary, when $\g$ is semisimple, the Gerstenhaber bracket
induces the null bracket on $H_c^*(\g)$. Applying the preceding
results, we give a complete description of the super Poisson bracket
in $\extgg$, and of the \gla $H_c^*(\g)$, when $\g = \gl(n)$
(\ref{6.2.6}).

We develop in Section \ref{Section7}, the theory of quadratic $2k$-Lie
algebra structures on a semisimple Lie algebra $\g$, in relation with
cyclic cochains (\ref{7.1.3}). This is a direct generalization of the
case of quadratic Lie algebras studied in Section \ref{Section4}. We
show that any invariant even cyclic cochain $F$ defines a quadratic
$2k$-Lie algebra (\ref{7.2.1}) and that $\extgg = H^\star (\g)$ is
contained in $H^\star(F)$. Finally, we give an interpretation of some
interesting examples given in \cite{Per} of $2k$-Lie algebras in terms
of the techniques developed in the present paper, pointing out where
these examples come from. Finally, we give some examples of quadratic
$2k$-Lie algebra structures on $\gl(n)$ (\ref{8.3.1}).

\subsection*{Acknowledgments} We wish to tank Jim Stasheff for
his valuable and enlightening suggestions that allowed us to greatly
improve the first version of this paper.

\section{$\Mcc(\g)$, $\Mcc_a(\g)$ and $2k$-Lie algebra
  structures}\label{Section1}

This Section is essentially a review, except \ref{1.5}. For more
details, see \cite{Gerst} and \cite{Nij-Rich}.

Let $\g$ be a complex vector space. We denote by $\Mcc(\g)$ the space
of multilinear mappings from $\g$ to $\g$. The space $\Mcc(\g)$ is
graded as follows:
\[ \Mcc(\g) = \sum_{k \geq 0} \Mcc^k (\g)\]
where $\Mcc^0 (\g) = \g$, $\Mcc^k (\g) = \{ F \colon \g^k \rightarrow
\g \mid F \ k\text{-linear} \}$, for $k \geq 1$.

\subsection{} \label{1.2}

The theory of associative algebra structures on $\g$ can be described
in a graded Lie algebra framework \cite{Gerst, Nij-Rich}: first,
consider $\Mcc(\g)$ with shifted grading $\Mcc^k[1] = \Mcc^{k+1} (\g)$
and denote it $\Mcc[1]$. Then define a graded Lie bracket on $\Mcc[1]$
as follows: for all $F \in \Mcc^p [1]$, $G \in \Mcc^q[1]$, then $[F,
  G] \in \Mcc^{p+q} [1]$ with
\begin{eqnarray*}
[F,G] (X_1, \dots, X_{p+q+1}) : = & \\ (-1)^{pq} \sum_{j=1}^{p+1}
(-1)^{q(j-1)} & F(X_1, \dots, X_{j-1}, G(X_j, \dots, X_{j+q}),
X_{j+q+1}, \dots, X_{p+q+1}) \\ - \sum_{j=1}^{q+1} (-1)^{p(j-1)} &
G(X_1, \dots, X_{j-1}, F(X_j, \dots, X_{j+p}), X_{j+p+1}, \dots,
X_{p+q+1}).
\end{eqnarray*}
for $X_1, \dots, X_{p+q+1} \in \g$.

Notice that when $X \in \Mcc^{-1} [1] = \g$, then $[X,G]$ is defined
by:
\begin{equation*}
[X,G] (X_1, \dots, X_q ) = - \sum_{j=1}^{q+1} (-1)^{j-1} \ G(X_1,
\dots, X_{j-1}, X, X_j, \dots, X_q).
\end{equation*}

Notice also that when $F$ and $G$ are in $\Mcc^0 [1] = \End ( \g)$,
then $[F, G]$ is the usual bracket of the two linear maps $F$ and $G$.

Now, suppose that $F \in \Mcc^1 [1]$ defines a product on $\g$ by:
\[ X \cdot Y = F(X,Y), \forall \  X, Y \in \g. \]

This product is associative if and only if $[F, F] = 0$. In this case,
the derivation $\ad( F)$ of the graded Lie algebra $\Mcc [1]$
satisfies $(\ad (F))^2 = 0$, so it defines a complex on $\Mcc (\g)$
which turns out to be the Hoschild cohomology complex of the
associative algebra defined by $F$ \cite{Gerst}.

\subsection{} \label{1.3}

In the remaining of the paper, we use $\Sk_{p,q}$ to denote the set of
all $(p,q)$-{\em unshuffles}, that is, elements $\sigma$ in the
permutation group $\Sk_{p+q}$ satisfying $\sigma(1) < \dots <
\sigma(p)$ and $\sigma(p+1) < \dots < \sigma (p+q)$.

The theory of Lie algebra structures on $\g$ can also be described in
a graded Lie algebra framework \cite{Gerst, Nij-Rich}. First, let
$\Mcc_a = \Mcc_a (\g)$ be the space of skew symmetric elements in
$\Mcc (\g)$. One has $\Mcc_a = \dss \sum_{k \geq 0} \Mcc_a^k$ with
$\Mcc_a^0 = \g$ and $\Mcc_a^1= \End (\g)$. Then consider $\Mcc_a$ with
shifted grading denoted by $\Mcc_a [1]$, and define a graded Lie
bracket as follows: for all $F \in \Mcc^p_a [1]$, $G \in \Mcc^q_a[1]$,
then $[F, G]_a \in \Mcc^{p+q}_a [1]$ with
\begin{eqnarray*}
 [F,G]_a \ (X_1, \dots, X_{p+q+1}) : = & \\ (-1)^{pq} \underset{\sigma
   \in \Sk_{q+1,p}}{\sum} & \e (\sigma) \ F(G(X_{\sigma(1)}, \dots,
 X_{\sigma(q+1)}), X_{\sigma(q+2)}, \dots, X_{\sigma(p+q+1)})\\ -
 \underset{\sigma \in \Sk_{p+1,q}}{\sum} & \e (\sigma)
 \ G(F(X_{\sigma(1)}, \dots, X_{\sigma(p+1)}), X_{\sigma(p+2)}, \dots,
 X_{\sigma(p+q+1)})
\end{eqnarray*}
for $X_1, \dots, X_{p+q+1} \in \g$.

Notice that when $X \in \Mcc^{-1}_a [1] = \g$, then
\begin{equation*}
[X,G]_a \ (X_1, \dots, X_q ) = - G(X, X_1, \dots, X_q) \; \; (= -
\io_X (G) (X_1, \dots, X_q)).
\end{equation*}

Moreover, when $F$, $G \in \Mcc_a^0 [1] = \End (\g)$, then $[F, G]_a$ is the
usual bracket of the linear maps $F$ and $G$.

Now, any $F \in\Mcc_a^1 [1]$ defines a bracket on $\g$ by
\[ [X, Y] = F(X,Y), \forall \  X, Y \in \g. \]

The Jacobi identity is satisfied if and only if $[F, F]_a = 0$.  In
this case, the derivation $\ad(F)$ of the graded Lie algebra $\Mcc_a
[1]$ satisfies $(\ad(F))^2 = 0$, so it defines a complex on $\Mcc_a$
which turns out to be the Chevalley cohomology complex with
coefficients in the adjoint representation, of the Lie algebra
structure defined by $F$.

At this point, let us quickly explain the relations between the two
brackets defined in \ref{1.2} and \ref{1.3}. First, define the {\em
  skew symmetrization map} $\Ao \colon \Mcc (\g) \rightarrow \Mcc_a
(\g)$:
\[ \Ao (F) (X_1, \dots, X_k) = \sum_{\sigma \in \Sk_k} \e (\sigma)
\ F(X_{\sigma(1)}, \dots, X_{\sigma(k)})\] with $F \in \Mcc^k(\g)$ and
$X_1, \dots, X_k \in \g$. One has:

\begin{prop}\label{1.4.2}
For all $F$, $G \in \Mcc(\g)$, $\Ao ([F,G]) = [ \Ao(F), \Ao (G)]_a$.
\end{prop}

Obviously, when $F \in \Mcc^1 [1]$ induces an associative product on
$\g$, then $\Ao (F)$ induces a Lie algebra structure on $\g$. However
one should notice that from Proposition \ref{1.4.2}, Lie algebra
structures of type $\Ao (F)$ can be obtained from a ``product'' $F$ on
$\g$ satisfying other conditions than associativity, for instance:

\begin{prop}
Let $F \in \Mcc^1 [1]$ such that there exists $\tau \in \Sk_3$
satisfying \linebreak $\tau . [F,F] = - \e (\tau) \ [F,F]$. Then $\Ao
(F)$ defines a Lie algebra structure on $\g$.
\end{prop}

\subsection{} \label{1.5} 
Let us introduce a concept of generalized Lie algebra structures on
$\g$: 

\begin{defn}
An element $F \in \Mcc^{2k-1}_a [1]$ is a {\em $2k$-Lie algebra
  structure on $\g$} if
\[[F, F]_a = 0.\]
\end{defn}

We shall often use a bracket notation: for $X_1, \dots, X_{2k} \in
\g$,
\[[X_1, \dots, X_{2k}] = F (X_1, \dots,X_{2k}).\]

The identity $[F, F]_a = 0$ can be seen as a generalized Jacobi
identity (see \cite{Filip, Per}).

Given a $2k$-Lie algebra structure $F$ on $\g$, $\ad (F)$ is an odd
derivation of $\Mcc_a [1]$ and satisfies $(\ad (F))^2 = 0$, so there
is an associated cohomology defined by \linebreak $\ker (\ad (F)) /$ $
\im (\ad (F))$, which can be interpreted as a generalization of the
Chevalley complex of \ref{1.3}.

\section{$\Dc(\g)$, $\Wc(\g)$, cohomology of $2k$-Lie algebras and
  $\GPB$}\label{Section2}

In this Section, with exception made to \ref{2.3} and \ref{3.2}, we
recall classical material needed in the paper.

\subsection{} \label{2.1}

We denote by $\Dc = \Dc(\g)$ the space of (graded) derivations of
$\extg$. The space $\Dc$ is graded by $\Dc = \sum_{k=-1}^n \Dc^k$ with
$D \in \Dc^d$ if $D (\Wedge^p \g) \subset \Wedge^{p+d} \g$, for all
$p$, and has a graded Lie algebra structure with the bracket defined
by:
\begin{equation}\label{2.1.3}
[D, D'] = D \circ D' - (-1)^{dd'} D' \circ D, \forall \ D \in \Dc^d,
D' \in \Dc^{d'}.
\end{equation}

We denote by $\iota_X$, $X \in \g$, the elements of $\Dc^{-1}$ defined
by 
\[\iota_X (\Omega) (Y_1, \dots, Y_k) := \Omega (X, Y_1, \dots, Y_k),
\forall \ \Omega \in \Wedge {}^{k+1} \g, X, Y_1, \dots, Y_k \in \g \ (k
\geq 0),\]

and $\iota_X(1) = 0$. When $\g$ is finite dimensional, given a basis
$\{X_1, \dots, X_n \}$ and its dual basis $\{ \omega_1, \dots,
\omega_n \}$, any element $D \in \Dc$ can be written in a unique way:
\[ D = \sum_{r=1}^n D_r \wedge \iota_{X_r} \]
where $D_r = D (\omega_r)$. Moreover, $\Dc$ is a simple Lie
superalgebra (often denoted by $W(n)$, see \cite{Scheu}) and there
exists an obvious vector space isomorphism $\Dt \colon \Mcc_a [1]
\rightarrow \Dc$ defined as $\Dt (\Omega \otimes X) = - \Omega \wedge
\iota_X$, $\forall \ \Omega \in \extg$, $X \in \g$ which turns out to
be a \gla isomorphism. 

Since we do not want to restrict ourselves to the finite dimensional
case, we give a proof of the following result:

\begin{prop}\label{2.1.4}
There exists a one to one \gla homomorphism $\Dt \colon \Mcc_a
[1] \rightarrow \Dc$ such that 
\[ \Dt (\Omega \otimes X) = - \Omega \wedge \iota_X, \forall
\ \Omega \in \extg, X \in \g.\]

When $\g$ is finite dimensional, $\Dt$ is an isomorphism.
\end{prop}

\begin{proof}
Given a basis $\{ X_r \mid r \in \Rc \}$ of $\g$, and the forms
$\omega_r$, $r \in \Rc$, defined by $\omega_r (X_s) = \delta_{rs}$,
$\forall r,s$, for $F \in \Mcc_a^k$, let $\Dt (F) = - \sum_{r \in \Rc}
{}^t F(\omega_r) \wedge \iota_{X_r}$. It is easy to see that though
its indexes set is infinite, this sum applied to an element $\Omega
\in \Wedge^w \g$ gives:
\begin{eqnarray*}
\Dt (F) (\Omega) (Y_1, \dots, Y_{k+w-1}) =& \\ - \sum_{\sigma \in
  \Sk_{k, w-1} } & \e(\sigma) \Omega(F(Y_{\sigma(1)}, \dots,
Y_{\sigma(k)}), Y_{\sigma(k+1)}, \dots, Y_{\sigma (k+w-1)}),
\end{eqnarray*}

for all $Y_1, \dots, Y_{k+w-1} \in \g$. It results that our definition
of $\Dt$ does not depend on the basis of $\g$, and that $\Dt (A
\otimes X) = - A \wedge \iota_X$, $A \in \Wedge^k \g$, $X \in
\g$. Keeping in mind the remark about the sum defining $\Dt$, we
compute for $G \in \Mcc_a^{k'}$:
\begin{eqnarray*}
&&[\Dt(F), \Dt(G)] = \\ &&\sum_{r,s}  \left( {}^t F(\omega_r) \wedge
\iota_{X_r} ({}^tG( \omega_s)) - (-1)^{(k+1)(k'+1)} {}^t G(\omega_r)
\wedge \iota_{X_r} ({}^tF( \omega_s)) \right) \wedge \iota_{X_s}.
\end{eqnarray*}

By a direct computation:
\begin{eqnarray*}
&\sum_{r} \left( {}^t F(\omega_r) \wedge {}^tG( \omega_s)) -
  (-1)^{(k+1)(k'+1)} {}^t G(\omega_r) \wedge \iota_{X_r} ({}^tF(
  \omega_s)) \right) (Y_1, \dots, Y_{k+k'-1}) \\ &= \omega_s \left(
  \sum_{\sigma \in \Sk_{k, k'-1} } \e(\sigma) G(F(Y_{\sigma(1)},
  \dots, Y_{\sigma(k)}), Y_{\sigma(k+1)}, \dots, Y_{\sigma (k+k'-1)})
  \right. - \\ & \left.  (-1)^{(k+1)(k'+1)} \sum_{\sigma \in \Sk_{k',
      k-1} } \e(\sigma) F(G(Y_{\sigma(1)}, \dots, Y_{\sigma(k')}),
  Y_{\sigma(k'+1)}, \dots, Y_{\sigma (k+k'-1)}) \right) \\ &= -{}^t
  [F,G](\omega_s)(Y_1, \dots, Y_{k+k'-1}),
\end{eqnarray*}

Hence:
\begin{eqnarray*}
[\Dt(F), \Dt(G)] = - \sum_{s} \ {}^t [F,G] (\omega_s) \wedge
\iota_{X_s} = \Dt([F,G])
\end{eqnarray*}

\end{proof}

In the sequel, given $F \in \Mcc_a (\g)$, we denote by $\Dt_F$ the
associated derivation of $\extg$. If $\g$ is finite dimensional, for
$D \in \Dc$, we denote by $\Ft_D$ the associated element in $\Mcc_a
(\g)$. Here are some examples:


\begin{ex}
If $T \in \End (\g) = \Mcc^0_a [1]$, then
\[\Dt_T (\Omega) (Y_1, \dots, Y_p) = - \sum_{i=1}^p \Omega (Y_1,
\dots, Y_{i-1}, T(Y_i), Y_{i+1}, \dots, Y_p)\] for all $\Omega \in
\Wedge^p \g$, $Y_1, \dots, Y_p \in \g$.
\end{ex}

\begin{ex}
If $F \in \Mcc^1_a [1]$, then
\begin{eqnarray}
\label{2.1.7}  \Dt_F (\Omega) (Y_1, &\dots&, Y_{p+1}) = \\
\notag && \sum_{i < j} (-1)^{i+j} \ \Omega (F(Y_i, Y_j), Y_1, \dots,
\hat{Y}_i, \dots, \hat{Y}_j, \dots, Y_{p+1})
\end{eqnarray}

for all $\Omega \in \Wedge^p \g$, $Y_1, \dots, Y_{p+1} \in \g$.
\end{ex}

\subsection{} \label{2.2}

Let $F$ be a Lie algebra structure on $\g$, then $F \in \Mcc_a^1 [1]$
and $[F, F]_a = 0$. Let $\partial = \Dt_F$, then $[\partial, \partial]
= 0$ gives $\partial^2 = 0$ and formula (\ref{2.1.7}) shows that the
associated complex in the Grassmann algebra $\extg$ is exactly the
Chevalley cohomology complex of trivial cohomology of $\g$. One
defines $\theta_X$:
\begin{equation*}\label{2.2.1}
\theta_X = [\iota_X, \partial] = \Dt_{\ad (X)}.
\end{equation*}

If $\{ X_r \mid r \in \Rc \}$ is a basis of $\g$, consider the forms
$\omega_r$, $r \in \Rc$, defined by $\omega_r (X_s) = \delta_{rs}$,
$\forall r,s$. The map $\theta$ defines a Lie algebra representation
of $\g$ in $\extg$ and one has:
\begin{equation}\label{2.2.2}
\partial = \frac{1}{2} \sum_{r \in \Rc} \omega_r \wedge \theta_{X_r}.
\end{equation}

Let us precise that this formula is well-known when $\g$ is finite
dimensional (see \cite{Koszul1}), and that a proof in the infinite
dimensional case is given in the proof of Lemma \ref{s.3.7} of the
present paper. In any case, a very important consequence of formula
(\ref{2.2.2}) is that any invariant in $(\extg)^\g$ is a cocycle.

\subsection{}\label{2.3}

Let us now check how \ref{2.2} can be extended to $2k$-Lie algebra
structures on $\g$. Let $F \in \Mcc_a^{2k-1} [1]$. Assume that $[F,
  F]_a = 0$, and let
\begin{equation}\label{2.3.1}
[Y_1, \dots, Y_{2k}] : = F (Y_1, \dots, Y_{2k})
\end{equation}
for $Y_1, \dots, Y_{2k} \in \g$. Denote by $ D = \Dt_F$ the associated
derivation of $\extg$. Using Proposition \ref{2.1.4}, one concludes
$D^2 = 0$, so one can define an associated cohomology $H^\star(F) =
\ker( D) / \im (D)$. One has
\begin{equation*}\label{2.3.2}
D \omega (Y_1, \dots, Y_{2k}) = - \omega ( [Y_1, \dots, Y_{2k}])
\end{equation*}
for $\omega \in \gd$, $Y_1, \dots, Y_{2k} \in\g$. We shall come back
to cohomology of $2k$-Lie algebras in Section \ref{Section7}.

In the remaining of this Section, we will assume that $\g$ is a finite
dimensional vector space with $\dim \g = n$. We state next some
properties of the Schouten bracket. For more details, we refer to
\cite{Koszul2}.

\subsection{}\label{3.1}

Let $\Wc = \Wc(\g) = \Pc \otimes \extg $, graded by $\Wc^p = \Pc
\otimes \Wedge^p \g$, where $\Pc$ is the symmetric algebra of $\gd$ .
Elements of $\Wc$ act as skew symmetric multivectors on $\Pc$ as
follows: for $\Omega \in \Wedge^p \g$, $P \in \Pc$, $f_1, \dots, f_p
\in \Pc$,
\[ (P \otimes \Omega) (f_1, \dots, f_p)_\varphi = P(\varphi) \ \Omega
((df_1)_\varphi, \dots, (df_p)_\varphi). \]

For instance, if $\{ X_1, \dots, X_p \}$ is a basis of $\g$ and
$\{\omega_1, \dots, \omega_p \}$ its dual basis, one has for all $i =
1, \dots, p$:
\[ \omega_i(f) = \dfrac{\partial f}{\partial X_i}, \forall f \in \Pc.\]

There is a natural $\wedge$-product on $\Wc$, defined by: for all $P,
P' \in \Pc$, $\Omega, \Omega' \in \extg$:
\[ (P \otimes \Omega) \wedge (P' \otimes \Omega') = P P' \otimes (\Omega
\wedge \Omega').\]

Each $f \in \Pc$ defines a derivation $\iota_f$ of degree $-1$ of $\Wc$ by:
\begin{eqnarray*}
\iota_f(P \otimes \Omega)(f_1, \dots, f_{p-1})_\varphi &=& P(\varphi)
\ \iota_{(df)_\varphi} (\Omega) ((df_1)_\varphi, \dots,
(df_{p-1})_\varphi)\\ & = & P(\varphi) \ \Omega((df)_\varphi,
(df_1)_\varphi, \dots, (df_{p-1})_\varphi).
\end{eqnarray*}

For instance, if $V \in \Wc^1$, one has $\iota_f(V) = V(f)$. There is
a graded Lie bracket on $\Wc[1]$ called the {\em Schouten bracket},
and defined by: for all $W \in \Wc^p [1]$, $W' \in \Wc^q[1]$, then
$[W, W']_S \in \Wc^{p+q} [1]$ with
\begin{eqnarray*}
  [W,W']_S \  (f_1, \dots, f_{p+q+1})  = & \\
  (-1)^{pq} \underset{\sigma \in \Sk_{q+1,p}}{\sum} & \e (\sigma) \ 
  W(W'(f_{\sigma(1)}, \dots, f_{\sigma(q+1)}), f_{\sigma(q+2)}, \dots,
  f_{\sigma(p+q+1)})\\ - \underset{\sigma \in \Sk_{p+1,q}}{\sum}
  & \e (\sigma) \ W'(W(f_{\sigma(1)}, \dots, f_{\sigma(p+1)}), f_{\sigma(p+2)},
  \dots, f_{\sigma(p+q+1)})
\end{eqnarray*}
for $f_1, \dots, f_{p+q+1} \in \Pc$.

Then for all $P, P' \in \Pc$, $\Omega \in \Wedge^{p+1} \g$, $\Omega'
\in \Wedge^{q+1} \g$:
\[ [ P \otimes \Omega, P' \otimes \Omega']_S = (-1)^{pq} P \otimes (\Omega'
\wedge \iota_{P'}(\Omega)) - P' \otimes (\Omega \wedge \iota_P (\Omega')).\]

As a particular case, one has $[\Omega, \Omega']_S = 0$, for all $\Omega$,
$\Omega' \in \extg$.

Let $W \in \Wc^1[1]$, then $W$ defines a Poisson bracket on $\Pc$ by
$\{P, P' \} = W(P,P')$ if and only if $[W, W]_S =0$. More generally,
as proposed in \cite{Per}, one can define {\em Generalized Poisson
  Brackets} ($\GPB$) as follows:

\begin{defn}\label{3.1.1}
An element $W \in \Wc^{2k-1} [1]$ is a $\GPB$ if $[W, W]_S = 0$.
\end{defn}

(see \cite{Per} where these structures are introduced and applications
are proposed).

\subsection{}\label{3.2}

Let us now show that $2k$-Lie algebras have associated $\GPB$, exactly
as Lie algebras have associated Poisson brackets. This will be a
consequence of the following construction: define a map $V \colon \Dc
= \Dc(\g) \rightarrow \Wc$ by $V_D = V(D) := -X \otimes \Omega$ for $D
= \Omega \wedge \iota_X$ with $\Omega \in \extg$, $X \in \g$. Then, it
is easy to check that:

\begin{prop}\label{3.2.2}
One has $V_{[D, D']} = [V_D, V_{D'}]_S$, for $D$, $D' \in
\Dc$. Moreover $V$ is a one to one graded Lie algebras homomorphism
from $\Dc$ into $\Wc[1]$.
\end{prop}

For example, given a $2k$-Lie algebra structure $F$ on $\g$, denoted
by $[Y_1, \dots, Y_{2k}]$ $ = F(Y_1, $ $\dots, Y_{2k})$, $\forall
\ Y_1, \dots, Y_{2k} \in \g$, let $D$ be the associated derivation
(see Proposition \ref{2.1.4}) in $\Dc$. Then one has:
\[V_D(f_1, \dots, f_{2k})_\varphi = \langle \varphi | [(df_1)_\varphi, \dots,
  (df_{2k})_\varphi] \rangle,\]

and since $[F,F]_a = 0$ by (\ref{2.1.4}), one has $[D, D] =0$ .  Using
Proposition \ref{3.2.2} above, $[V_D, V_D]_S = 0$, so $V_D$ defines a
$\GPB$ on $\Pc$.

Finally, using \ref{2.1} and Proposition \ref{3.2.2}, one deduces an
inclusion of the simple Lie superalgebra $W(n)$ into the graded Lie
algebra $\Wc[1]$, endowed with the Schouten bracket which provides a
natural realization of $W(n)$.

\section{Application to identities of standard polynomials, and
  cohomology}\label{Sections}

In this Section, $\g$ denotes an associative algebra, with product
$m$. We also use the notation: $X . Y = m (X,Y)$, $\forall \ X,Y \in
\g$. We assume that $m$ has a unit $\mathbf{1}_m$, but this is not
really necessary.

\subsection{}\label{s.1}

We first define the iterated $m_k$ ($k \geq 0$) of $m$ as:
\[ m_0 = \mathbf{1}_m, \ m_1 = \Id_\g, \ m_2 = m, \dots, \ m_k (Y_1,
\dots, Y_k) = Y_1 . \ldots . Y_k, \ \forall \ Y_1, \dots, Y_k \in \g,
\dots \]

It is easy to check that:

\begin{prop}\label{s.1.2}
For all $k$, $k' \geq 0$, one has:
\begin{eqnarray*} 
\ [ m_{2k}, m_{2k'} ] & = & 0, \\ \ [ m_{2k}, m_{2k'+1} ] & = & \ (2 k
-1) \ m_{2 k+ 2k'} , \\ \ [ m_{2k+1}, m_{2k'+1} ] & =& \ 2 (k - k')
\ m_{2k + 2k' +1}.
\end{eqnarray*}
\end{prop}

Hence the space generated by $\{ m_k, k \geq 0\}$ is a subalgebra of
the \gla $\Mcc (\g)$ of Section \ref{Section1}.

\subsection{}\label{s.2}

Now define the {\em standard polynomials} $\Ac_k$ ($k \geq 0$) on $\g$
as:
\[ \Ac_k := \Ao(m_k)\]

Using Propositions \ref{1.4.2} and \ref{s.1.2}, one immediately
obtains:

\begin{prop}\label{s.2.2}
For all $k$, $k' \geq 0$, one has:
\begin{eqnarray*} 
\ [ \Ac_{2k}, \Ac_{2k'} ]_a & = & 0, \\ \ [ \Ac_{2k}, \Ac_{2k'+1} ]_a
& = & \ (2 k -1) \ \Ac_{2k+2k'} , \\ \ [ \Ac_{2k+1}, \Ac_{2k'+1} ]_a &
=& \ 2 (k - k') \ \Ac_{2k + 2k' +1}.
\end{eqnarray*}
\end{prop}

Let $\Ac$ be the subspace generated by $\{ \Ac_{k}, k \geq 0
\}$. Hence $\Ac$ is a subalgebra of the graded Lie algebra $\Mcc_a
(\g)$ of Section \ref{Section1}. The standard polynomial $\Ac_2$ is
the Lie algebra structure on $\g$ associated to $m$. Since $[
  \Ac_{2k}, \Ac_{2k} ]_a = 0$, $\forall k$, we conclude:

\begin{prop}\label{s.2.5}
The standard polynomials $\Ac_{2k}$, $k \geq 1$ define $2k$-Lie
algebra structures on $\g$.
\end{prop}

Remark that $\Ac_k$ is a $\g$-invariant map from $\g^k$ to $\g$ for
the Lie algebra structure. Moreover the standard polynomial $\Ac_{2k}$
is a coboundary of the adjoint representation of the Lie algebra $\g$
whereas $ [ \Ac_{2}, \Ac_{2k-1} ] = \Ac_{2k}$.

\subsection{}\label{s.3}

Let us now define an associative product on $\Ac$. First consider an
associative product $\circ$ on $\Mcc (\g)$:
\begin{eqnarray*}
(F \circ G) (Y_1, \dots, Y_{p+q}) = F(Y_1, \dots, Y_p) . G(Y_{p+1},
  \dots, Y_{p+q}), 
\end{eqnarray*} 
for all $F \in \Mcc^p (\g), G \in \Mcc^q(\g)$, $Y_1, \dots, Y_{p+q}
\in \g$.

Then define an associative product $\times$ on $\Mcc_a (\g)$ by:
\[ (F \times G)(Y_1, \dots, Y_{p+q}) = \underset{\sigma \in
    \Sk_{p,q}}{\sum} \e (\sigma) \ F(Y_{\sigma(1)}, \dots,
Y_{\sigma(p)}) \ . \ G(Y_{\sigma(p+1)}, \dots, Y_{\sigma(p+q)}),\] for
all $F \in \Mcc_a^p (\g), G \in \Mcc_a^q(\g)$, $Y_1, \dots, Y_{p+q}
\in \g$. By a straightforward computation, one has:

\begin{prop}\label{s.3.3}
For all $F$, $G \in \Mcc_a(\g)$, $\Ao (F \circ G) = \Ao(F) \times
\Ao(G)$.
\end{prop}

It is obvious that $m_k = \underbrace{m_1 \circ \dots \circ m_1}_{k
  \textrm{ times}}$, so:

\begin{cor}\label{s.3.4}
$\Ac_k = \underbrace{\Ac_1 \times \dots \times \Ac_1}_{k \textrm{
      times}}$, for all $k \geq 1$ and $\Ac_k \times \Ac_{\ell} =
  \Ac_\ell \times \Ac_k = \Ac_{k + \ell}$, for all $k, \ell \geq 0$.
\end{cor}

As a consequence, $\Ac$ is a commutative algebra for the
$\times$-product.

Any element $Z \in \g$ defines a super derivation $\iota_Z$ of degree
$-1$ of the $\times$-product of $\Mcc_a(\g)$ by: for all $F \in
\Mcc_a^p (\g)$, $Y_1, \dots, Y_{p-1} \in \g$,
\[ \iota_Z (F) (Y_1, \dots, Y_{p-1}) := F(Z, Y_1, \dots, Y_{p-1}),\]

Denote by $\Zo(\g)$ the center of the algebra $\g$. If $Z \in \Zo
(\g)$, one has $\iota_Z (\Ac_2) = 0$. Hence using Corollary
\ref{s.3.4} and the derivation property of $\iota_Z$, we deduce:

\begin{prop}\label{s.4.2}
Assume that $Z \in \Zo (\g)$. Then for all $k$, 
\[ \iota_Z (\Ac_{2k} ) = 0 \text{   and   } \iota_Z (\Ac_{2k+1}) = Z \ . \
\Ac_{2k}.\]
\end{prop}

This Proposition expresses classical identities on standard
polynomials, generally written in the case $Z = \mathbf{1}_m$.

Let us now assume that $\g$ is equipped with a trace, that is, a
linear form $\Tr \colon \g \rightarrow \CC$ satisfying:
\[ \Tr(X . Y) = \Tr(Y.X), \ \forall \ X, Y \in \g. \]

Let $\Wedge \g$ be the Grassmann algebra of $\g$. We extend the trace
$\Tr$ to a map $\Tr \colon \Mcc_a (\g) \rightarrow \Wedge \g$ defined
by:
\[ \Tr (F) (Y_1, \dots, Y_p) =  \Tr (F (Y_1, \dots, Y_p)), \]
for all $F \in \Mcc_a^p (\g)$, $Y_1, \dots, Y_p \in \g$.

\begin{prop}
One has $\Tr (F \times G) = (-1)^{pq} \ \Tr (G \times F)$, for all $F \in
\Mcc_a^p(\g)$, $G \in \Mcc_a^q(\g)$.
\end{prop}

\begin{proof}
Let $F \in \Mcc_a^p (\g), \ G \in \Mcc_a^q(\g), \ Y_1, \dots, Y_{p+q}
\in \g$:
\begin{eqnarray*}
\Tr (F \times G)(Y_1, \dots, Y_{p+q}) = &\\ = \underset{\sigma \in
  \Sk_{p,q}}{\sum} \e (\sigma) \ & \Tr(F(Y_{\sigma(1)}, \dots,
Y_{\sigma(p)}) \ . \ G(Y_{\sigma(p+1)}, \dots, Y_{\sigma(p+q)})) \\ =
\underset{\sigma \in \Sk_{p,q}}{\sum} \e (\sigma) \ & \Tr(
G(Y_{\sigma(p+1)}, \dots, Y_{\sigma(p+q)}) \ .  \ F(Y_{\sigma(1)},
\dots, Y_{\sigma(p)}) )
\end{eqnarray*}

Given $\sigma \in \Sk_{p,q}$, define $\tau \in \Sk_{q,p}$ as $\tau(1)
= \sigma(p+1)$, $\dots$, $\tau(q) = \sigma(p+q)$ and $\tau(q+1) =
\sigma(1)$, $\dots$, $\tau(q+p) = \sigma(p)$. Then one has $\e (\tau)
= (-1)^{pq} \e(\sigma)$, so:
\begin{eqnarray*}
&& \Tr (F \times G)(Y_1, \dots, Y_{p+q}) = \\ &=& (-1)^{pq}
  \underset{\tau \in \Sk_{q,p}}{\sum} \e (\tau) \ \Tr( G(Y_{\tau(1)},
  \dots, Y_{\tau(q)}) \ . \ F(Y_{\tau(q+1)}, \dots, Y_{\tau(p+q)})
  )\\ &=& (-1)^{pq} \Tr (G \times F) (Y_1, \dots, Y_{p+q}).
\end{eqnarray*}

\end{proof}

Hence our extension of the trace has, in fact, the properties of a
$\extg$-valued super trace on the graded algebra $(\Mcc_a (\g),
\times)$. Denoting the super bracket associated to the
$\times$-product on $\Mcc_a (\g)$ by:
\[[F, G]_\times = F \times G - (-1)^{pq} G \times F, \forall  \ F \in \Mcc_a^p (\g),G
\in \Mcc_a^q (\g),\] 

one obtains
\begin{cor}  \label{s.3.5}
$\Tr ([F, G]_\times ) = 0$.
\end{cor}

\begin{prop} \label{s.3.6}
One has $\Tr (\Ac_{2k}) = 0$ ($k \geq 1$) and $\Tr (\Ac_{2k+1})$ ($k
\geq 0$) is an invariant cocycle for the (trivial) cohomology of the
Lie algebra $\g$.
\end{prop}

\begin{proof}
For the first claim, use $[\Ac_1, \Ac_{2k -1} ]_\times = 2 \ \Ac_{2k}$
and apply the Corollary above. For the second, we remark that
$\Ac_{2k+1}$ is a $\g$-invariant map from $\g^{2k+1}$ into $\g$, so
$\Tr(\Ac_{2k+1}) \in (\extg)^\g$ and therefore a cocycle by the
following classical Lemma.
\end{proof}

\begin{lem}  \label{s.3.7}
Let $\hk$ be a Lie algebra. Then any invariant cochain in $(\Wedge
\hk)^\hk$ is a cocycle.
\end{lem}

\begin{proof}
If $\hk$ is finite dimensional, the result is well known
(\cite{Koszul1}) and is a direct consequence of the formula $\partial
= \frac12 \sum_{i=1}^n \omega_i \wedge \theta_{X_i}$ where $\partial$
is the differential, $\{X_1, \dots, X_n \}$ a basis of $\hk$ and $\{
\omega_1, \dots, \omega_n \}$ its dual basis.

For the sake of completeness, we give a proof in the general case, let
$\{X_i \mid i \in I \}$ be a basis of $\hk$, and $\{ \omega_i \mid i
\in I \}$ be the forms defined by $\omega_i(X_j) = \delta_{ij}$,
$\forall \ i, j$. We claim that the formula $\partial = \frac12
\sum_{i \in I} \omega_i \wedge \theta_{X_i}$ is still valid. To prove
this, let $D = \frac12 \sum_{i \in I} \omega_i \wedge
\theta_{X_i}$. Though its indexes set is infinite, this sum exists
since for $\Omega \in \Wedge^p \hk$ and $Y_1, \dots, Y_{p+1} \in \hk$,
one has:
\[ \frac12 \sum_{i \in I} \omega_i \wedge \theta_{X_i} (\Omega) (Y_1, \dots,
Y_{p+1}) = \frac12 \sum_{j=1}^{p+1} (-1)^{j+1} \sum_{i \in I}
\omega_i(Y_j) \theta_{X_i}(\Omega) (Y_1, \dots, \hat{Y_j}, \dots,
Y_{p+1}) \]
Then
\begin{eqnarray*}
\notag D(\Omega) (Y_1, \dots, Y_{p+1}) = - \frac12 \sum_{j=1}^{p+1}
(-1)^{j+1} \scriptstyle{\times}
\end{eqnarray*}
\begin{eqnarray*}
\notag \left( \sum_{k=1}^{j-1} \Omega (Y_1, \dots, [Y_j, Y_k], \dots,
\hat{Y_j}, \dots, Y_{p+1}) + \sum_{k=j+1}^{p+1} \Omega (Y_1, \dots,
\hat{Y_j}, \dots, [Y_j, Y_k], \dots, Y_{p+1}) \right)
\end{eqnarray*}
\begin{eqnarray*}
\notag &=& \frac12 \sum_{j=1}^{p+1} (-1)^j \left( \sum_{k < j}
(-1)^{k+1} \Omega ([Y_j, Y_k], Y_1, \dots, \hat{Y_k}, \dots,
\hat{Y_j}, \dots, Y_{p+1}) + \right. \\ \notag && \phantom{@@@@@@@@@}
\left. \sum_{j<k} (-1)^k \Omega ([Y_j, Y_k], Y_1, \dots, \hat{Y_j},
\dots, \hat{Y_k}, \dots, Y_{p+1}) \right) \\ \notag &=& \sum_{j < k}
(-1)^{j+k} \Omega ([Y_j, Y_k], Y_1, \dots, \hat{Y_j}, \dots,
\hat{Y_k}, \dots, Y_{p+1}) = \partial( \Omega) (Y_1, \dots, Y_{p+1})
\end{eqnarray*}
\end{proof}

Now recall the well-known formula (e.g. \cite{Kostant1}):

\begin{prop}\label{s.5.7}
$\Tr(\Ac_{2k+1} (Y_1, \dots, Y_{2k+1})) = (2k+1) \Tr(\Ac_{2k} (Y_1,
  \dots, Y_{2k}). Y_{2k+1})$, for all $Y_1, \dots, Y_{2k+1} \in \g$.
\end{prop}

This formula will be reinterpreted in Section \ref{Section5} in terms
of cyclic cohomology of the Lie algebra $\g$: $\Ac_{2k}$ is a cocycle
of the adjoint action (actually a coboundary since $[\Ac_2,
  \Ac_{2k-1}] = \Ac_{2k}$), and Proposition \ref{s.5.7} tells that it
is a cyclic cocycle, as will be defined in Section \ref{Section5}.

\begin{ex}\label{s.6}
Assume that $\g = \gl(n)$. Then $H^\star(\g)$ can be completely
described in terms of standard polynomials (see e.g. \cite{Kostant1}
or \cite{Fuchs}):
\[ H^\star (\g) = \ext [\Tr(\Ac_1), \Tr(\Ac_3), \dots, \Tr(\Ac_{2n-1})]. \]

Moreover, by the Amitsur-Levitzki theorem (\cite{Amitsur, Kostant1}):
\[ \Ac_k = 0, \textrm{  if  } k \geq 2n.\]

So $\dim \Ac = 2n$. For the $\times$-product, $\Ac \simeq \CC[X] /
X^{2n}$. For the graded bracket, the structure of $\Ac$ is given by
formula (\ref{2.3.1}), $\Ac_2$ is the Lie algebra structure on $\g$
and the standard polynomials $\Ac_4, \dots, \Ac_{2n-2}$ define
$2k$-Lie algebra structures on $\g$ by Proposition \ref{s.2.2}.
\end{ex}

\begin{ex}\label{s.7}
More generally, let $V$ be an infinite dimensional vector space. Let
$\g$ be the space of finite rank linear maps. So $\g$ is an ideal of
the associative algebra $\End(V)$. There is a vector spaces
isomorphism $\g \simeq V^* \otimes V$ defined by $(\omega \otimes v)
(v') = \omega(v') \ v$, for all $v$, $v' \in V$, $\omega \in V^*$. So
we can define the trace $\Tr(X)$ when $X \in \g$ by $\Tr(\omega
\otimes v) := \omega(v)$, for all $\omega \in V^*$, $v \in V$. It is
easy to check that $\Tr([X,Y]) = 0$, for all $X$, $Y \in \g$, so the
preceding results apply. Moreover, the symmetric bilinear form $B$
defined on $\g$ by $B(X,Y) = \Tr(XY)$ is non degenerate and invariant,
therefore $\g$ is a quadratic Lie algebra. Since $\gl(n) \subset \g$,
$\forall \ n$, by Example \ref{s.6}, we can conclude that
\[ \ext[\Tr(\Ac_1), \Tr(\Ac_3), \dots, \Tr(\Ac_{2n-1}), \dots ]
\subset H^\star (\g)\]

\begin{prop}\label{s.8}
Let $a_{2n+1}=\Tr(\Ac_{2n+1})$. Then 
\[H^\star(\g) = \ext[a_1, a_3, \dots,  a_{2n+1}, \dots ]. \]
\end{prop}

\begin{proof}
Recall that for any Lie algebra $\hk$, there is an isomorphism
$H^k(\hk) \simeq H_k(\hk)^*$, induced by the restriction $\Omega \in
Z^k(\hk) \mapsto \Omega|_{Z_k(\hk)}$ where $H_k(\hk) $ is the homology
of $\hk$ defined as $H_k(\hk) = Z_k(\hk) / B_k(\hk)$ (with $Z_k(\hk)$
the cycles and $B_k(\hk)$ the boundaries).

Let us define $\Sc = \{ S = (W, W') \mid W, W'$ complementary
subspaces of $ V$ with $\dim(W) < \infty \}$ and for $S = (W, W') \in
\Sc$, $\g_S = \{ X \in \g \mid X(V) \subset W, \ X(W') = \{0 \}
\}$. Then $\g_S$ is a subalgebra of the (associative or Lie) algebra
$\g$ and one has $\g_S \simeq \gl(\dim(W))$. It is easy to check that
given $X_1, $, $\dots$, $X_r \in \g$, there exists $S \in \Sc$ such
that $X_i \in \g_S$, $\forall i = 1, \dots, r$. It results that, if $c
\in \ext^k(\g)$, there exists $S$ such that $c \in \ext^k(\g_S)$, so
that $\ext^k(\g) = \cup_{S \in \Sc} \ext^k(\g_S)$.

Set $\Ec = \ext[a_1, a_3, \dots, a_{2n+1}, \dots] \subset H^\star(\g)$
and $\Ec^k = \Ec \cap H^k (\g)$. Then $\dim(\Ec^k)$ $= \sharp \ I_k$
with $I_k = \{ (i_j) \in \{0,1\}^\NN \mid \sum_{j \in \NN} (2 j +1)
\ i_j = k \}$. We fix a basis $\{ \Omega_i \mid i \in I_k \}$ of
$\Ec^k$.

Given $c \in Z_k(\g)$, denote by $\overline{c}$ its class in
$H_k(\g)$. Let us assume that $\Omega_i(\overline{c}) = 0$, $\forall i
\in I_k$. Take $S \in \Sc$ such that $c \in \ext^k(\g_S)$, then by
(\ref{s.6}), $\{ \Omega_i \mid i \in I_k \}$ generates $H^k(\g_S) =
H_k(\g_S)^*$ and since $c \in Z_k(\g_S)$, it results that $c \in
B_k(\g_S) \subset B_k(\g)$, therefore $\overline{c} = 0$. So, $\{
\Omega_i \mid i \in I_k \}$ is free in $H_k(\g)^*$ and $\cap_{i \in
  I_k} \ker(\Omega_i) = \{0\}$. It results that $\dim(H_k(\g)) =
\sharp \ I_k$. Since $H^k(\g) = H_k(\g)^*$, one has $\dim(H^k(\g)) =
\sharp \ I_k$ and since $\Ec^k \subset H^k(\g)$, one obtains $\Ec^k =
H^k(\g)$.
\end{proof}

\begin{rem} \label{s.9}
From $H^1(\g) = \CC \ \Tr$, we deduce that $[\g, \g] =
\ker(\Tr)$. From $H^2(\g) = \{ 0 \}$, we deduce that $\g$ has no (non
trivial) central extension.
\end{rem}

\end{ex}

\section{Super Poisson brackets and quadratic Lie algebras} \label{Section4}

The canonical Poisson bracket on $\RR^{2n}$ appears as the leading
term of a quantization of the algebra of polynomial functions by the
Weyl algebra, the so-called {\em Moyal product}. We will develop a
similar formalism, replacing polynomials (i.e. commuting variables) by
skew multilinear forms (i.e. skew commuting variables) and the Weyl
algebra by the Clifford algebra. The leading term of the deformation
will be the {\em super Poisson bracket}.

\subsection{} \label{4.1}

Let us give a definition of the Clifford algebra that is well adapted
to the realization of this algebra as a deformation of the exterior
algebra. Denote by $\Cc_t$, $t \in \CC$, the associative algebra with
basis $\{ e_I, I \in \ZZ_2^n \}$ and product defined by
\[ e_I \star e_J = (-1)^{\Omega(I,J)} \ t^{|IJ|} \ e_{I \underset{(2)}{+} J} \]
where $\Omega$ is the bilinear form associated to the matrix
$(a_{ij})_{i,j=1}^n$ with $a_{ij} = 1$ if $i>j$ and $0$ otherwise.

Take $I_i=(j_k) \in \ZZ_2^n$, with $j_i=1$ and $0$ otherwise. Set
$e_i= e_{I_i}$, $i=1, \dots, n$ and $V = \spa \{e_1, \dots, e_n
\}$. When $t=0$, one obtains $\Cc_0 = \ext(V)$. When $t \neq 0$,
$\Cc_t$ is the {\em Clifford algebra}. The following relations hold:
\[e_i^2 = t, \ \forall \ i, \quad e_i \star e_j + e_j \star e_i = 0, \ i \neq j,
\] 
\[ e_{i_1} \star e_{i_2} \star \dots \star e_{i_p} = e_{i_1} \wedge e_{i_2}
\wedge \dots \wedge e_{i_p}, \textrm{ if } i_1 < i_2 < \dots < i_p\]

So that $\Cc_t$ is the quotient algebra of the tensor algebra $T(V)$
by the relations:
\[ v \otimes v = t . \ B(v,v) . \ 1, v \in V, \] where $B$ is the
bilinear form $B(e_i, e_j) = \delta_{ij}$, for all $i$, $j$, and we
recover the usual definition of the Clifford algebra.

But, we are mainly interested in realizing $\Cc_t$ as a deformation of
$\ext(V)$. Using:
\[ t^k = 0^k + t \ \delta_{k,1} + t^2 \ \delta_{k,2} + \dots \ . \]

this deformation becomes transparent:

\begin{prop}\label{4.1.6}
One has \[ e_I \star e_J = e_I \wedge e_J + \sum_{k=1}^n t^k
\ D_k(e_I, e_J)\] where $D_k(e_I, e_J)$ $=$ $\delta_{|IJ|,k}
\ (-1)^{\Omega(I,J)} \ e_{I \underset{(2)}{+} J}$.
\end{prop}

Symmetry properties of the coefficients are resumed in:

\begin{prop}
For all $\Omega \in \ext^w(V)$, $\Omega' \in \ext^{w'}(V)$,
\[D_j (\Omega, \Omega') = (-1)^j (-1)^{ww'} D_j(\Omega', \Omega).\]
\end{prop}

We insist on the fact that $\Cc_t$ is not a $\ZZ$-graded, but only a
$\ZZ_2$-graded algebra. The associated Lie superalgebra has bracket:
\[ [\Omega, \Omega']_\star = 2 \ \sum_{p\geq 0} \ t^{2p+1} \ D_{2p+1} (\Omega,
\Omega').\]

\begin{defn} \label{4.1.3}
We define the {\em super Poisson bracket} on $\ext(V)$ by:
\[ \{ \Omega, \Omega' \} = 2 \ D_1 (\Omega, \Omega'), \ \forall
\ \Omega, \Omega' \in \ext(V).\]
\end{defn}

Since $[., .]_\star$ satisfies the super Jacobi identity, so does
$\{., .\}$. Moreover, since $\ad_\star( \Omega)$ is derivation of the
$\Cc_t$-product, $\ad_{\tt{P}} (\Omega) := \{ \Omega, .\}$ is a
derivation of the $\wedge$-product (actually of degree $(w -2)$ if
$\Omega \in \ext^w(V)$).

Finally, by a straightforward computation, one gets:
\begin{eqnarray} 
\{ v_1 \wedge \dots \wedge v_p, w_1 \wedge \dots &\wedge& w_q \} = 2
\ (-1)^{p+1} \scriptstyle{\times} \label{4.1.11} \\ \underset{j=1,
  \dots, q}{\underset{i=1, \dots, p}{\sum}} (-1)^{i+j} \ B(v_i, w_j) &
v_1& \wedge \dots \wedge \hat{v_i} \wedge \dots \wedge v_p \wedge w_1
\wedge \dots \wedge \hat{w_j} \wedge \dots \wedge w_q, \notag
\end{eqnarray}
for all $v_1, \dots, v_p$, $w_1, \dots, w_q \in V$.

Comparing with the formulas given in \cite{Scheu}, we conclude that
the Lie superalgebra $\ext(V) / \CC$ is isomorphic to the simple Lie
superalgebra $\tilde{H} (n)$. Notice that $\ext(V) / \CC \simeq
\ad_{\tt{P}} (\ext(V)) \subset \Der (\ext(V)) = \Dc(V^*)$, so we
obtain the classical inclusion $\tilde{H}(n) \subset W(n)$
(\cite{Scheu}).

\subsection{} \label{4.2}

Let us modify slightly the formalism in \ref{4.1} in order to apply it
to Lie algebras deformation theory. We begin with a $n$-dimensional
vector space $\g$ and we set $V = \gd$. We assume that $\g$ is a
quadratic space with bilinear form $B$. Denote by $\{X_1, \dots, X_n
\}$ an orthonormal basis of $\g$ and by $\{ \omega_1, \dots, \omega_n
\}$ the dual basis; we define $B$ on $\gd$ by $B(\omega_i, \omega_j) =
\delta_{ij}$. Applying the construction in \ref{4.1} with $e_i =
\omega_i$, $i = 1, \dots, n$, we get a super Poisson bracket on
$\extg$ and it is easy to check that:

\begin{prop} \label{4.2.1}
For all $\Omega \in \Wedge^w \g$, $\Omega' \in \extg$, one has 
\[ \{ \Omega, \Omega'\} = 2 \ (-1)^{w+1} \ \sum_{j=1}^n \iota_{X_j} (\Omega)
\wedge \iota_{X_j} (\Omega').\]
\end{prop}

This formula is valid in any orthonormal basis of $\g$ and it is
enough for our purpose in Section \ref{Section4}, but a general
formula can be found in Lemma \ref{4bis.5.2}. There is a Moyal type
formula which gives the Clifford product in terms of the super Poisson
bracket: let $m_{\scriptscriptstyle{\wedge}}$ be the product from
$\extg \otimes \extg \to \extg$, and define $\Fc \colon \extg \otimes
\extg \to \extg \otimes \extg$ by
\[  \Fc(\Omega \otimes  \Omega') =  \ (-1)^w \ \sum_{j=1}^n
  \iota_{X_j} (\Omega) \otimes \iota_{X_j} (\Omega')\] for all $\Omega
  \in \Wedge^w \g$, $\Omega' \in \extg$. Then:

\begin{prop}\label{4.2.2} \hfill
\[\Omega \star \Omega' = m_{\scriptscriptstyle{\wedge}} \circ \exp(-t
  \Fc)(\Omega \otimes \Omega').\]
\end{prop}

\begin{proof}
As in the beginning of \ref{4.2}, let $e_i = \omega_i$ and let
$\partial_i = \iota_{X_i}$, $i = 1, \dots, n$. As in \ref{4.1}, for $I
= (i_1, \dots, i_n) \in \ZZ_2^n$, let $e_I = e_1^{i_1} \wedge \dots
\wedge e_n^{i_n}$ and $\partial_I = \partial_1^{i_1} \circ \dots \circ
\partial_n^{i_n}$. For $J = (j_1, \dots, j_n) \in \ZZ_2^n$, let $J^I =
j_1^{i_1} \dots j_n^{i_n}$. One has $\partial_I(e_J) = (-1)^{\Omega
  (I,J)} J^I e_{I \underset{(2)}{+} J}$. Since all $\partial_i \otimes
\partial _i$ commute, and $\partial_i^2 = 0$, one has:
\[ \left( \sum_i \partial_i \otimes \partial_i \right)^k = k! \sum_{|I| = k}
\partial_I \otimes \partial_I\]

For $k > 0$, one has:
\begin{eqnarray*}
&&m_\wedge \circ \Fc^k(e_R \otimes e_S) \\ & =& (-1)^{|R|}
  (-1)^{|R|-1} \dots (-1)^{|R|-(k-1)} m_\wedge \circ \left( \sum_i
  \partial_i \otimes \partial_i \right)^k (e_R \otimes e_S) \\& =&
  (-1)^{k|R|} (-1)^{\frac{k(k-1)}{2}} k! \sum_{|I|= k}
  (-1)^{\Omega(I,R)} (-1)^{\Omega(I,S)} R^I S^I \ e_{I
    \underset{(2)}{+} R} \wedge e_{I \underset{(2)}{+} S}
\end{eqnarray*}

This vanishes, except if $k = |RS|$, and in that case, the only
remaining term in the sum is when $I = RS$. We compute this term:
\[m_\wedge \circ \Fc^k(e_R \otimes e_S) = (-1)^{k|R|}
  (-1)^{\frac{k(k-1)}{2}} k! (-1)^{\Omega(RS,R)+\Omega(RS,S)}
(-1)^{\Omega(RS+R,RS+S)} e_{R \underset{(2)}{+} S}\]

But one has $\Omega(A,B) + \Omega(B,A) = |A||B|- |AB|$, so:
\[\Omega(RS,R)+\Omega(RS,S) = k|R|-k, \quad \text{ and } \quad
\Omega(RS,RS) = \frac{k(k-1)}{2}.\]

So finally, we have proved that 
\[m_\wedge \circ \Fc^k(e_R \otimes e_S) = (-1)^{k} k!
(-1)^{\Omega(R,S)} \delta_{|RS|,k} e_{R \underset{(2)}{+} S}\]

On the other hand, by Proposition \ref{4.1.6}, one has
\[ e_R \star e_S = e_R \wedge e_S + (-1)^{\Omega(R,S)} \sum_{k=1}^n
\delta_{|RS|,k} e_{R \underset{(2)}{+} S}\]

so the result follows.                                            
\end{proof}

\begin{rem}
An equivalent formula is given in \cite{Kostant2}, but without the use
of super Poisson bracket.
\end{rem}

\subsection{} \label{4.3}

A derivation $D \in \Dc$ is {\em Hamiltonian} if it belongs to
$\ad_{\tt{P}}(\extg)$. Actually, the space of Hamiltonian derivations
is a subalgebra of $\Dc$, that we denote by $\Hc(\g)$, which is
isomorphic to $\extgq = \Wedge\g / \CC$ and therefore, by
(\ref{4.1.11}), isomorphic to the simple Lie superalgebra
$\widetilde{H}(n)$. Here a simple characterization of Hamiltonian
derivations:

\begin{prop} \label{4.3.1}
A derivation $D = \sum_r D_r \wedge \iota_{X_r}$ is Hamiltonian if and
only if $\iota_{X_r} (D_s) + \iota_{X_s} (D_r) = 0 $, $\forall \ r,s$.
\end{prop}

\begin{proof}
When the condition is satisfied, one has $D = \ad_{\tt{P}}(\Omega)$
where $\Omega = \frac{1}{2 w} \sum_r D_r \wedge \omega_r$ and $w =
\deg(D) + 2$.
\end{proof}

\begin{rem}
A Hamiltonian derivation is a derivation of the $\wedge$-product and
also of the super Poisson bracket. 
\end{rem}

In fact, one has:

\begin{prop}
Let $D \in \Dc$. Then $D$ is Hamiltonian if and only if $D$ is a
derivation of the super Poisson bracket.
\end{prop}

\begin{proof}
Let $D = \sum_r D_r \wedge \iota_{X_r}$ with $D \in \Wedge^d \g$, then
$D_r = D(\omega_r)$. Since $\{ \omega_r, \omega_s \} \in \CC$,
assuming that $D$ is a derivation of the super Poisson bracket, one
has:
\[ 0 = D( \{ \omega_r, \omega_s \}) = 2 \ (-1)^{d+1} \left(
\iota_{X_r} (D_s) + \iota_{X_s} (D_r) \right) \] and the result
follows by Proposition \ref{4.3.1}.
\end{proof}

\subsection{} \label{4.4-0}

We now want to apply super Poisson brackets to the theory of quadratic
Lie algebras, in a deformation framework that we will quickly set
up. Given a quadratic Lie algebra $(\g_0, B_0)$ with bilinear form
$B_0$ and product $[.,.]$, a {\em deformation} $(\g_t, B_t)$ of
$(\g_0, B_0)$ is:

\begin{enumerate}

\item a deformation $\g_t$ of $\g$ in the usual sense, so:
\[ [X,Y]_t = [X,Y] + t C_1(X,Y) + \dots, \forall \ X, Y \in \g,\]

\item a formal bilinear form $B_t = B_0 + t B_1 + \dots$ such that
\[ B_t([X,Y]_t,Z) = -B_t(Y, [X,Z]_t), \forall \ X, Y, Z \in \g,\]

\end{enumerate}

Two deformations $(\g_t, B_t)$ and $(\g'_t, B'_t)$ with respective
brackets $[.,.]_t$ and $[.,.]'_t$ are {\em equivalent} if there exists
$T_t = \Id + t T_1 + \dots$ such that:
\[ [X,Y]'_t = T_t^{-1} ([T_t(X), T_t(Y)]) \quad \text{ and } \quad
B'_t(X,Y) = B_t(T_t(X), T_t(Y)), \forall \ X, Y \in \g.\]

\begin{prop}\label{4.4.1}
Any deformation $(\g_t, B_t)$ of $(\g_0, B_0)$ is equivalent to a
deformation with unchanged bilinear form.
\end{prop}

\begin{proof}
Fix an orthonormal basis $\{ e_1, \dots, e_n \}$ of $\g$ with respect
to $B_0$. By a Gram-Schmidt type strategy, one can construct $\{
e_1(t), \dots, e_n(t) \}$ such that:
\[ e_\ell (t)  = \lambda_1 (t) e_1(t) + \dots + \lambda_{\ell - 1} (t)
e_{\ell -1}(t)+ e_\ell, \forall \ell \leq n,\]

with $\lambda_j(t) \in t \ \CC[[t]]$, and $B_t(e_\ell (t) , e_m(t))=
0$, for all $\ell, m \leq n$. Since \linebreak $[B_t(e_\ell(t),
  e_\ell(t))]_{t=0} = B_0(e_\ell, e_\ell) = 1$, $\forall \ell \leq n$,
$B_t(e_\ell(t), e_\ell(t))$ is invertible, and
\[e_\ell' (t) = \frac{1}{\left( B_t(e_\ell (t) ,
  e_\ell(t))\right)^{\frac12}} \ e_\ell(t)\] 
does satisfy $B_t(e'_\ell (t) , e'_m(t)) = \delta_{\ell m }$, $\forall
\ell, m$.

Now if we define $T_t$ by $T_t(e_\ell) = e_\ell'(t)$, $\forall \ell
\leq n$, and a new deformation
\[ [X,Y]_t  = T_t([T_t(X), T_t(Y)]_t), \ \forall X, Y \in \g,\]

with bilinear form $B'_t(X,Y) = B_t(T_t(X), T_t(Y)) = B_0 (X,Y)$,
$\forall X$, $Y \in \g$, we obtain a deformation that is equivalent to
the initial one.
\end{proof}

So if one want to study quadratic Lie algebras in terms of deformation
theory, one can restrict to quadratic Lie algebras with a specified
bilinear form, and that is what we shall do next.

\subsection{} \label{4.4}

The constructions made in the beginning of this section can now be
applied as follows: given a finite dimensional quadratic Lie algebra
$\g$ with bilinear form $B$, let $\partial $ be the corresponding
derivation of $\extg$ (i.e. the differential of the trivial cohomology
complex of $\g$, see \ref{2.2}), we define:
\[ I(X,Y,Z) := B([X,Y], Z), \ \forall \ X, Y, Z \in \g\]

Then one has:

\begin{prop} \label{4.4.2}
\hfill

\begin{enumerate}

\item $I \in (\Wedge^3 \g)^\g$.

\item $\partial = - \frac{1}{2} \ad_{\tt{P}}(I)$.

\item $\{ I, I \} = 0$.

\end{enumerate}

\end{prop}

\begin{proof} 
The assertion (1) is obvious. To show (2), let $\{X_1, \dots, X_n \}$
be an orthonormal basis of $\g$ and $\{ \omega_1, \dots, \omega_n \}$
the dual basis. Then for all $Y$, $Z \in \g$:
\begin{eqnarray*}
-\frac12 \adpo (I) (\omega_i) (Y,Z) = - \left( \sum_j \iota_{X_j} (I)
\wedge \iota_{X_j} (\omega_i) \right) (Y,Z) = - B([X_i, Y],Z) =
\\ \notag = - B(X_i, [Y,Z]) = - \omega_i ([Y,Z]) = \partial \omega_i
(Y,Z)
\end{eqnarray*}

Hence, $\partial = - \frac12 \adpo (I)$.

Finally $\adpo(\{I,I\}) = [ \adpo(I), \adpo(I)] = 4 [ \partial,
  \partial] = 8 \partial^2 = 0$. So $\{I, I \} = 0$ and that proves
(3).
\end{proof}

Note that $\partial$, $\iota_X$ and $\theta_X= [\iota_X, \partial]$,
$\forall \ X \in \g$ are all Hamiltonian derivations.

\subsection{} \label{4.5}

Conversely, assume that $\g$ is (only) a finite dimensional quadratic
vector space. Fix $I \in \Wedge^3 \g$ and define $\partial = -
\frac{1}{2} \ad_{\tt{P}} (I)$. Then the formula
\[ \ad_{\tt{P}} \left( \{ \Omega, \Omega' \} \right) = [ \ad_{\tt{P}}
  (\Omega), \ad_{\tt{P}} (\Omega') ], \forall \ \Omega, \Omega' \in
\extg\]

leads to
\begin{equation} \label{4.5.2} 
[\partial,\partial]= 0 \text{ if and only if } \{ I, I \}=0.
\end{equation}

Let $F =\Ft_\partial$ be the structure on $\g$ associated to
$\partial$ (see \ref{2.1} and \ref{2.2}), then from (\ref{4.5.2}), it
follows:

\begin{prop} \label{4.5.3}
$F$ is a Lie algebra structure if and only if $\{I, I \} = 0$. In that
  case, with the notation $[X, Y] = F(X,Y)$, one has:
\[I(X, Y, Z)= B([X, Y], Z), \ \forall \ X, Y, Z \in \g, \]

the form $B$ is invariant and $\g$ is a quadratic Lie algebra.
\end{prop}

\begin{proof} 
We have to prove that if $F$ is a Lie algebra structure, then $I(X, Y,
Z)= B([X, Y], Z)$, $\forall \ X, Y, Z \in \g$.

Let $\{X_1, \dots, X_n \}$ be an orthonormal basis, then $\partial = -
\sum_k \iota_{X_k}(I) \wedge \iota_{X_k}$, so $F = \sum_k
\iota_{X_k}(I) \otimes X_k$, and therefore $B([X_i, X_j], X_k) =
\iota_{X_k}(I) (X_i, X_j) = I(X_i, X_j, X_k)$, for all \linebreak $i$,
$j$, $k$.
\end{proof}

\begin{rem}\label{4.6.2}
Using \ref{4.4-0}, \ref{4.4} and \ref{4.5}, it appears that $\extg[2]$
with super Poisson bracket is a \gla associated to deformation theory
of finite dimensional quadratic Lie algebras: by \ref{4.4-0}, one can
assume that $B$ does not change, then quadratic Lie algebra structures
with the same $B$ are in one to one correspondence with elements $I
\in \Wedge^3 \g$ such that $\{I, I\} = 0$ (\ref{4.4}, \ref{4.5}). An
equivalent description can be given in terms of Hamiltonian
derivations, i.e. of the \gla $\Hc(\g) = \adpo(\extg) \simeq \extg /
\CC = \extgq$.

Let us note that in this picture, one has to redefine equivalence: a
priori, one might think that equivalence should be defined as Lie
algebras isomorphism keeping $B$ fixed. But this is too restrictive,
since $[.,.]$ and $\lambda(t) [.,.]$, with $\lambda (t) = 1 + t
(\dots)$ will not be equivalent in that sense as they should be. So
one has rather to work with the notion of a conformal equivalence,
i.e.  an equivalence defined by a Lie algebras isomorphism $T(t) = \Id
+ t (\dots)$satisfying $B(T_t(X), T_t(Y)) = \mu(t) B(X,Y)$, with $\mu
(t) = 1 + t (\dots)$. This will change the corresponding \gla: one can
consider the subalgebra $\CC R \oplus \Hc(\g)$ of $\Dc(\g)$ (where $R
= \sum_i \omega_i \wedge \iota_{X_i}$ is the {\em super radial vector
  field}), rather than $\Hc(\g)$. Hence, there are some adaptations to
carry out, which will not be developed here since they are somewhat
standard. Let us only indicate that in this framework, if $(\g_0,B_0)$
is the initial quadratic Lie algebra with associated $I_0 \in
\Wedge^3(\g)$, then the first obstruction to triviality of a quadratic
deformation will lie in $H^3(\g) / \CC \ I_0$. For instance, if $\g_0$
is semisimple, it is shown in \cite{Koszul1} that $H^3(\g_0)$ and the
space of symmetric invariant bilinear forms on $\g_0$ are isomorphic,
the isomorphism being $B \mapsto I_B$ where $I_B(X,Y,Z) = B([X,Y],Z)$,
$\forall \ X$, $Y$, $Z \in \g$. It results that when $\g_0$ is simple,
it is rigid in quadratic deformation theory.
\end{rem}

\section{Elementary quadratic Lie algebras}\label{Section4bis}

\subsection{}

First, let us recall two results:

\begin{prop}\label{4bis.1}
Let $V$ be a finite dimensional vector space and $I$ a $k$-form in
$\Wedge^k V$. Denote by $V_I$ the orthogonal subspace in $V^*$ of the
subspace $\{ X \in V \mid \iota_X(I) = 0 \}$. Then $\dim(V_I) \geq k$
and if $I$ is non zero, $I$ is decomposable if and only if $\dim V_I=
k$. In this case, if $\{\omega_1, \dots, \omega_k \}$ is a basis of
$V_I$, one has $I = \alpha \ \omega_1 \wedge \dots \wedge \omega_k$,
for some $\alpha \in \CC$ (\cite{Bour}).
\end{prop}

\begin{prop}\label{4bis.1.2}
Let $V$ be a finite dimensional quadratic vector space with a non
degenerate symmetric bilinear form $B$. For a subspace $W$ of $V$,
denote by $W^\perp$ its orthogonal subspace in $V$ with respect to $B$
and $W^{\perp^*}$ its orthogonal in $V^*$. Let $\phi$ be the
isomorphism from $V$ onto $V^*$ induced by $B$.  Then
$\phi|_{W^\perp}$ is an isomorphism from $W^\perp$ onto $W^{\perp^*}$,
so $\dim(W^\perp) = \dim(V) - \dim(W)$. One has $V = W \oplus W^\perp$
if and only if $W \cap W^\perp = \{0 \}$ and in this case the
restriction of $B$ to $W$ or $W^\perp$ is non degenerate.
\end{prop}

In the rest of this Section, $\g$ will denote a finite dimensional
quadratic Lie algebra with bilinear form $B$. Denote by $\Zo(\g)$ the
center of $\g$ and by $I_\g$ the element of $\Wedge^3 \g$ defined by
$I_\g (X,Y,Z) = B([X,Y],Z)$, $\forall \ X$, $Y$, $Z \in \g$.

\begin{defn}\label{4.1.3ele}
We say that $\g$ is an {\em elementary} quadratic Lie algebra if
$I_\g$ is decomposable.
\end{defn}

Remark that the obvious identity $\Zo(\g)^\perp = [\g, \g]$ holds
here. As a consequence,

\begin{prop}\label{4bis.2.2}
Let $\g$ be a non Abelian quadratic Lie algebra. Then $\dim([\g, \g])
\geq 3$. Moreover, $\g$ is elementary if and only if the equality
holds.
\end{prop}

\begin{proof}
Since $\Zo(\g)^\perp = [\g, \g]$ and $V_{I_\g} = \Zo(\g)^{\perp^*}$,
the result follows directly from Proposition \ref{4bis.1}.
\end{proof}

\begin{cor}\label{4bis.2.3}
Let $\g$ be an elementary quadratic Lie algebra. Then all codjoint
orbits have dimension at most 2.
\end{cor}

\begin{proof}
Let $\omega \in \gd$ and $X_\omega \in \g$ such that $\omega =
\phi(X_\omega)$. Then $\ad(\g)(\omega) = \phi([\g,X_\omega]) \subset
\phi([\g,\g])$, so $\dim(\ad(\g)(\omega)) \leq 3$ and since all
coadjoint orbits have even dimension, the result is proved.
\end{proof}

\begin{rem} \label{4bis.2.4}
Suppose that $\g$ is (only) a finite dimensional quadratic vector
space and let $I$ be a decomposable $3$-vector in $\Wedge^3 \g$. Then
it is easy to check that $\{I, I \} = 0$ for the super Poisson
bracket. So by Proposition \ref{4.5.3} there is a quadratic elementary
Lie algebra structure on $\g$ such that $I(X,Y,Z) = B([X,Y],Z)$,
$\forall \ X$, $Y$, $Z \in \g$.
\end{rem}

In the sequel, we will classify all elementary non Abelian quadratic
Lie algebras. Our proof is a two steps one: first, in \ref{4bis.4}, we
demonstrate a result on quadratic Lie algebras that reduces the
classification problem to small dimensions, namely between 3 and
6. Then in \ref{4bis.5}, we proceed by classifying these small
dimensional elementary quadratic algebras. Explicit commutations in a
canonical basis with respect to $B$ are computed as well.

\subsection{} \label{4bis.4}

Here is the reduction result on quadratic Lie algebras:

\begin{prop}\label{4bis.4.1}
Let $\g$ be a non Abelian quadratic Lie algebra with bilinear form
$B$. Then there exist a central ideal $\zk$ and an ideal $\lk \neq \{
0 \}$ such that:

\begin{enumerate}
\item One has $\g = \zk \oplus \lk$, and $\lk$ and $\zk$ are
  orthogonal with respect to $B$.

\item The ideals $\zk$ and $\lk$ are quadratic (with bilinear forms
  induced by the restriction of $B$) and $\lk$ is non
  Abelian. Moreover, $\lk$ is elementary if and only if $\g$ is
  elementary.

\item The center $\Zo(\lk)$ is totally isotropic and one has
\[ \dim(\Zo(\lk)) \leq \frac12 \dim(\lk) \leq \dim([\lk, \lk]). \]
\end{enumerate}

\end{prop}

\begin{proof}
Let $\zk_0 = \Zo(\g) \cap [\g, \g]$. Fix any subspace $\zk$ such that
$\Zo(\g) = \zk_0 \oplus \zk$. Since $\Zo(\g)^\perp = [\g, \g]$, one
has $B(\zk_0, \zk) = \{ 0 \}$ and $\zk \cap \zk^\perp = \{0 \}$. It
results from Proposition \ref{4bis.1.2} that $\g = \zk \oplus \lk$
where $\lk = \zk^\perp$.

Since $B([\g,\g], \zk) = \{0\}$, one has $[\g,\g] \subset \lk$. It is
easy to check that $\Zo(\lk) = \zk_0$ and $[\lk, \lk] = [\g, \g] =
\Zo(\g)^\perp$ so $\Zo(\lk)$ is totally isotropic; moreover the
restriction of $B$ to $\zk$ is not degenerate, so $\lk$ is quadratic
and clearly non Abelian since $\zk$ is central in $\g$. If $\g$ is
elementary, $\lk$ is a fortiori elementary. If $\lk$ is elementary,
let $I_\lk = \omega_1 \wedge \omega_2 \wedge \omega_3$, $\omega_1$,
$\omega_2$, $\omega_3 \in \lk^*$, extend the $\omega_i$ to $\g$ by
$\omega_i|_\zk =0$. Since $I_\g (X,Y,Z) = 0$, $\forall X \in \zk$,
$Y$, $Z \in \g$, one concludes $I_\g = \omega_1 \wedge \omega_2 \wedge
\omega_3$, hence $\g$ is elementary. Finally $\Zo(\lk) \subset [\lk,
  \lk] = \Zo(\lk)^\perp$ implies $\dim(\lk) - \dim([\lk, \lk]) \leq
\dim([\lk, \lk])$ and the last inequality follows.
\end{proof}

\begin{cor}\label{4bis.4.2}
Let $\lk$ be an elementary non zero quadratic Lie algebra such that
$\Zo(\lk) $ is totally isotropic. Then one has
\[ 3 \leq \dim(\lk) \leq 6.\]
\end{cor}

\begin{proof}
Use Propositions \ref{4bis.4.1}(3) and \ref{4bis.2.2}.
\end{proof}

\subsection{} \label{4bis.5}

We shall now finish the classification of non Abelian elementary
quadratic Lie algebras. This classification is reduced, by Proposition
\ref{4bis.4.1} and Corollary \ref{4bis.4.2}, to the case of non zero  
elementary quadratic $\lk$ with a totally isotropic center
$\Zo(\lk)$. Applying Proposition \ref{4bis.4.2} one has $3 \leq
\dim(\lk) \leq 6$. Note that if $\dim(\lk) = 3$, one has $\lk =[\lk,
  \lk]$ (Proposition \ref{4bis.2.2}), so $\lk \simeq \slk(2)$ and $B$
is the Killing form up to a scalar. So we have to consider $\dim(\lk)
\geq 4$ (therefore $\dim(\Zo(\lk)) \geq 1$).

We need the following Lemma:

\begin{lem}\label{4bis.5.2}
Let $V$ be a quadratic vector space with bilinear form $B$. Define $B$
on $V^*$ by $B(\omega, \omega') := B(\phi^{-1} (\omega), \phi^{-1}
(\omega'))$, $\forall \ \omega$, $\omega' \in V^*$ ($\phi$ as in
Proposition \ref{4bis.1.2}). Let $\{\omega_1, \dots, \omega_n\}$ be a
basis of $V^*$, $\{X_1, \dots, X_n \}$ its dual basis and $\{Y_1,
\dots, Y_n \}$ the basis of $\g$ defined by $Y_i =
\phi^{-1}(\omega_i)$. Then the super Poisson bracket on $\extg$ is
given by
\[\{\Omega, \Omega'\}= 2 \ (-1)^{w+1} \sum_{i,j}
B(Y_i,Y_j) \iota_{X_i}(\Omega) \wedge \iota_{X_j}(\Omega'), \ \Omega
\in \Wedge{}^w \g, \Omega' \in \extg.\]
\end{lem}

\begin{proof}
Using Proposition \ref{4.2.1}, one has 
\[\{\Omega, \Omega'\}= 2 \ (-1)^{w+1} \sum_{i,j} \alpha_{ij}
\iota_{X_i}(\Omega) \wedge \iota_{X_j}(\Omega'),\]

$\Omega \in \Wedge^w \g$, $\Omega' \in \extg$ and $\alpha_{ij} =
\frac12 \{ \omega_i, \omega_j\}$. But from \ref{4.1}, one has
\linebreak $\{ \omega_i, \omega_j\} $ $=$ $ 2 B(\omega_i, \omega_j) $
$= 2 B(Y_i, Y_j)$.
\end{proof}

\begin{prop}\label{4bis.4.3}
Let $\lk$ be an elementary quadratic Lie algebra with non zero totally
isotropic center $\Zo(\lk)$. Then:

\begin{enumerate}
\item If $\dim(\lk) =6$, there exists a basis $\{Z_1, Z_2, Z_3, X_1,
  X_2, X_3 \}$ of $\lk$ such that:

\begin{itemize}
\item[(i)] $\{Z_1, Z_2, Z_3 \}$ is a basis of $\Zo(\lk)$.

\item[(ii)] $B(Z_i, Z_j) = B(X_i,X_j) = 0$, $B(Z_i, X_j) =
  \delta_{ij}$, $\forall \ i,j$.

\item[(iii)] $[X_1, X_2 ] = Z_3$, $[X_2, X_3 ] = Z_1$, $[X_3, X_1 ] =
  Z_2$ and the other brackets vanish.
   
\end{itemize}

\item If $\dim(\lk) =5$, there exists a basis $\{Z_1, Z_2, X_1, X_2, T
  \}$ of $\lk$ such that:

\begin{itemize}

\item[(i)] $\{Z_1, Z_2 \}$ is a basis of $\Zo(\lk)$.

\item[(ii)] $ B(Z_i, Z_j) = B(X_i, X_j) = 0$, $B(Z_i, X_j) =
  \delta_{ij}$, $\forall \ i,j$, $ B(T, Z_i) = B(T, X_i) = 0$, $ B(T,
  T) = 1$.

\item[(iii)] $[X_1, T] = - Z_2$, $[X_2, T ] = Z_1$, $[X_1, X_2 ] = T$
  and the other brackets vanish.
   
\end{itemize}

\item If $\dim(\lk) =4$, then $\dim(\Zo(\lk)) =1$ and there exist
  totally isotropic subspaces $\ik$ with basis $\{Z, P \}$ and $\ik'$
  with basis $\{X,Q \}$ such that $\Zo(\lk) \subset \ik \subset [\lk,
    \lk ]$, $\lk = \ik \oplus \ik'$ and~:

\begin{itemize}

\item[(i)] $\Zo(\lk) = \CC \ Z$, $B(Z,X) = B(P,Q) =1$, $B(Z,Q) =
  B(X,P) =0$.

\item[(ii)] $[X, P] = P$, $[X,Q] = -Q$, $[P,Q] = Z$ and the other
  brackets vanish.
\end{itemize}

\end{enumerate}

\end{prop}

\begin{proof}
\hfill

\begin{enumerate}

\item Assuming that $\dim(\lk) =6$, one has $\dim(\Zo(\lk)) = 3$, so
  $\Zo(\lk) = [\lk, \lk] = \Zo(\lk)^\perp$. Using \cite{Bour}, there
  is a totally isotropic subspace $\lk'$ such that $\lk = \Zo(\lk)
  \oplus \lk'$. With the notation of Proposition \ref{4bis.1.2}, since
  $\phi|_{\lk'}$ is an isomorphism from $\lk'$ onto $\Zo(\lk)^*$, we
  can find a basis $\{ Z_1, Z_2, Z_3 \}$ of $\Zo(\lk)$ and a basis $\{
  X_1, X_2, X_3 \}$ of $\lk'$ such that $B(Z_i,X_j) =
  \delta_{ij}$. Then
\[\Zo(\lk)^{\perp^*} = \spa\{ X_1^*, X_2^*, X_3^* \} = \spa \{
\phi(Z_1),\phi(Z_2), \phi(Z_3) \}.\]

Let $I_\lk = B([X,Y],Z)$, $\forall \ X, Y, Z \in \lk$. Since
$V_{I_\lk} = \Zo(\lk)^{\perp^*}$, it results from Proposition
\ref{4bis.1} that $I_\lk = \alpha \ X_1^* \wedge X_2^* \wedge X_3^*$,
$\alpha \in \CC$. Replacing $X_1 $ by $\frac1\alpha X_1$ and $Z_1$ by
$\alpha Z_1$, we can assume that $\alpha = 1$. Using Proposition
\ref{4.4.2} and Lemma \ref{4bis.5.2}, $\partial = - \frac1 2
\ad_{\tt{P}} (I) = - \sum_{i=1}^3 \iota_{X_i} (X_1^* \wedge X_2^*
\wedge X_3^*) \wedge \iota_{Z_i}$, so by \ref{2.2} and \ref{2.1},
$[X,Y] = \sum_{i=1}^3 \iota_{X_i} ( X_1^* \wedge X_2^* \wedge X_3^*)
(X,Y) \ Z_i$, $\forall \ X$, $ Y \in \lk$ and the commutation rules
follow.

\item Assuming $\dim(\lk) = 5$, one has $\dim(\Zo(\lk))=2$. Using
  \cite{Bour}, there is a totally isotropic subspace $\lk'$ and a
  one-dimensional subspace $\lk''$ such that $\lk = \Zo(\lk) \oplus
  \lk' \oplus \lk''$ and $B(\Zo(\lk) \oplus \lk', \lk'') = \{ 0
  \}$.Then one can find a basis $\{Z_1, Z_2 \}$ of $\Zo(\lk)$, a basis
  $\{X_1, X_2 \}$ of $\lk$ and a basis $\{T \}$ of $\lk''$ such that:
  $B(Z_i, X_j) = \delta_{ij}$, $\forall \ i$, $j$ and $B(T,T) =
  1$. Therefore 
\[\Zo(\lk)^{\perp^*} = \spa\{ X_1^*, X_2^*, T^* \}=\spa \{ \phi(Z_1),
\phi(Z_2), \phi(T) \}.\]

So $I_\lk = \alpha \ X_1^* \wedge X_2^* \wedge T^*$, $\alpha \in
\CC$. Replacing $X_1 $ by $\frac1\alpha X_1$ and $Z_1$ by $\alpha
Z_1$, we can assume that $\alpha = 1$. By Proposition \ref{4.4.2} and
Lemma \ref{4bis.5.2}, one obtains $\partial = - \frac1 2 \ad_{\tt{P}}
(I) = - \sum_{i=1}^2 \iota_{X_i} (X_1^* \wedge X_2^* \wedge T^*)
\wedge \iota_{Z_i}$ $- \iota_T (X_1^* \wedge X_2^* \wedge T^*) \wedge
\iota_T$, so by \ref{2.2} and \ref{2.1}, $[X,Y] = \sum_{i=1}^2
\iota_{X_i} ( X_1^* \wedge X_2^* \wedge T^*) (X,Y) \ Z_i + \iota_T
(X_1^* \wedge X_2^* \wedge T^*) \ T$, $\forall \ X$, $Y \in \lk$ and
the commutation rules follow.

\item Assuming $\dim(\lk) = 4$, one has $\dim(\Zo(\lk))=1$. Using
  \cite{Bour}, there is a totally isotropic 2-dimensional subspace
  $\ik$ such that $\Zo(\lk) \subset \ik$. Since $\Zo(\lk)^\perp =
       [\lk, \lk]$, one has $\ik \subset [\lk, \lk]$. Using
       \cite{Bour} once more, there exists a totally isotropic $\ik'$
       such that $\lk = \ik \oplus \ik'$. Let us write $\ik = \spa \{
       Z,P \}$, $\ik' = \spa \{ X,Q \}$ with $\Zo(\lk) = \CC \ Z$ and
       $B(Z,X) = B(P,Q) = 1$, $B(Z,Q) = B(X,P) =0$. Therefore
\[\Zo(\lk)^{\perp^*} = \spa\{ P^*, Q^*, X^* \} = \spa \{ \phi(Q),
\phi(P), \phi(Z) \}.\] 

So $I_\lk = \alpha \ P^* \wedge Q^* \wedge X^*$, $\alpha \in
\CC$. Replacing $P $ by $\frac1\alpha P$ and $Q$ by $\alpha Q$, we can
assume that $\alpha = 1$. Using Proposition \ref{4.4.2}, Lemma
\ref{4bis.5.2}, \ref{2.2} and \ref{2.1} as above, one finds $[A,B] = [
  \iota_P ( P^* \wedge Q^* \wedge X^*) \ Q + \iota_Q (P^* \wedge Q^*
  \wedge X^*) \ P + \iota_X (P^* \wedge Q^* \wedge X^*) \ Z] (A,B) $,
$\forall \ A$, $B \in \lk$ and the commutation rules follow.

\end{enumerate}

As a final remark, the brackets in (1), (2) and (3) do satisfy Jacobi
identity thanks to Remark \ref{4bis.2.4}.
\end{proof}

\begin{rem}
In the Proposition above, cases (1) and (2) are nilpotent Lie algebras
and case (3) is a solvable, non nilpotent Lie algebra, with derived
algebra the Heisenberg algebra.
\end{rem}

\section{Cyclic cochains and cohomology of quadratic Lie
  algebras}\label{Section5} 

\subsection{} 

First we need some notations: $\g$ will be a $n$-dimensional quadratic
vector space with bilinear form $B$ and $\extgp = \sum_{k \geq 1}
\Wedge^k \g$, which is an associative algebra without unit. If $\g$ is
a quadratic Lie algebra, we denote by $F_0$ its bracket (i.e. $F_0
(X,Y) = [X,Y]$, $X$, $Y \in \g$), by $\partial = \Dt (F_0)$ (see
\ref{2.2}) the differential of $\extg$, by $H^*(\g)$ the corresponding
cohomology, and by $H^*_+(\g)$ the restricted cohomology,
i.e. $H^*_+(\g)= \sum_{k \geq1} H^k(\g)$ which is an algebra without
unit (for the induced wedge product).

When $\g$ is a $n$-dimensional quadratic vector space, $\extg$ is a
\gla for the super Poisson bracket with grading $\extg [2]$. Denote by
$\extgq$ the quotient \gla $\extgq = \extg / \CC$, and by $[.,.]_Q$
its bracket. The map $\adpo \colon \extg \to \Dc(\g)$ is a \gla
homomorphism, we define the \gla $\Hc(\g)$ of Hamiltonian derivations
to be the image $\Hc(\g) = \adpo(\extg)$, as in \ref{4.3}. There is an
obvious \gla isomorphism from $\extgq$ onto $\Hc(\g)$, and since
$\extgq \simeq \Ht(n)$ (see \ref{4.3}), the \gla $\extgq$, $\Hc(\g)$
and $\Ht(n)$ are isomorphic. Moreover, if $\g$ is a quadratic Lie
algebra, since $\partial$ is Hamiltonian (see Proposition
\ref{4.4.2}), the super Poisson bracket induces a \gla structure on
$H^*(\g)$ and also on $H^*_Q(\g)= H^*(\g) / \CC$.

\subsection{} \label{5.1}

Given $C \in \Mcc_a^k (\g)$ (see \ref{1.3}), we define $\hat{C}$ by:
\begin{eqnarray*} 
& & \text{ if } k=0, \ C \in \g, \ \hat{C} (Y) : = B(C , Y), \forall
  \ Y \in \g, \label{5.1.1} \\ & & \text{ if } k > 0, \ \hat{C}(Y_1,
  \dots, Y_{k+1}) : = B(C(Y_1, \dots, Y_k) , Y_{k+1}), \ \forall
  \ Y_1, \dots, Y_{k+1} \in \g. \label{5.1.2}
\end{eqnarray*}

\begin{defn} \label{5.1.3}
$C$ is a {\em cyclic cochain} if
\[ \hat{C} (Y_1, \dots, Y_{k+1}) = (-1)^k \hat{C} (Y_{k+1}, Y_1,
\dots, Y_k), \forall \ Y_1, \dots, Y_{k+1} \in \g.\]
\end{defn}

We denote by $\Cc_c(\g)$ the space of cyclic cochains. 

\begin{prop} \label{5.1.4}

\hfill

\begin{enumerate}
\item $C$ is a cyclic cochain if and only if $\hat{C} \in \Wedge_+
  \g$. The map $\Theta$ from $\Cc_c(\g)$ into $\extgp$ defined by
  $\Theta(C) = \hat{C}$, is one to one.

\item When $\g$ is finite dimensional, the map $\Theta \colon
  \Cc_c(\g) \to \Wedge_+ \g$ is an isomorphism.

\item $\Cc_c(\g)$ is a subalgebra of the \gla $\Mcc_a(\g)$

\end{enumerate}

\end{prop}

\begin{proof}

\hfill

\begin{enumerate}
\item Let $\tau$ be the cycle $\tau = ( 1 \ 2 \ \dots \ k+1) \in
  \Sk_{k+1}$.  Given $\sigma \in \Sk_{k+1}$, let $\ell = \sigma^{-1}
  (k+1)$, then $\sigma' = \sigma \circ \tau^\ell \in \Sk_k$. If $C$ is
  cyclic, one has $\tau^{-1} . \hat{C} = \e(\tau) \ \hat{C}$. So
  $\sigma. \hat{C} = (\sigma' \circ \tau^{-\ell}). \hat{C} =
  \e(\sigma) \hat{C}$ and therefore $\hat{C} \in \Wedge_+ \g$. Since
  $B$ is non degenerate, $\Theta$ is clearly one to one.

\item Given $\Omega \in \Wedge^{k+1} \g$, define $D \in \Mcc_a^k (\g)$
  by $\Omega(Y_1, \dots, Y_k, Y) = B(D(Y_1, \dots, Y_k),Y)$, $\forall
  Y_1, \dots, Y_k,Y \in \g$. Then $\Omega = \hat{D}$.

\item Let $F \in \Mcc_a^p(\g)$, and $G \in \Mcc_a^q(\g)$, from (1) we
  have to prove that:
\[ B([F,G]_a (Y_1,  \dots, Y_{p+q-1}),Y_{p+q}) =  B([F,G]_a (Y_1,
\dots, Y_{p+q-2}, Y_{p+q}),Y_{p+q-1}),\] 

for all $Y_1, \dots, Y_{p+q}\in \g$. Using the formulas of \ref{1.2},
we can write the left hand side as a sum of four terms, $B([F,G]_a
(Y_1, \dots, Y_{p+q-1}),Y_{p+q}) = \alpha + \beta + \gamma + \delta$
where:
\begin{eqnarray*}
 \alpha = (-1)^{(p-1)(q-1)} \underset{\underset{\sigma(p+q-1) =
     p+q}{\sigma \in \Sk_{q,p-1}}}{\sum} (\dots) \quad &\text{ and }&
 \quad \beta = (-1)^{(p-1)(q-1)} \underset{\underset{\sigma(q) = p+q
     -1}{\sigma \in \Sk_{q,p-1}}}{\sum} (\dots)\\ \gamma = -
 \underset{\underset{\sigma(p+q-1) = p+q -1}{\sigma \in
     \Sk_{p,q-1}}}{\sum} (\dots) \quad &\text{ and }& \quad \delta = -
 \underset{\underset{\sigma(p) = p+q -1}{\sigma \in
     \Sk_{p,q-1}}}{\sum} (\dots)
\end{eqnarray*}

In $\alpha$, we can commute, up to a sign, $Y_{p+q-1}$ and
$Y_{p+q}$. In $\delta$, we commute, up to a sign, $F(Y_{\sigma(1)},
\dots, Y_{\sigma(p-1)}, Y_{p+q-1})$ and $Y_{p+q}$ to obtain:
\begin{eqnarray*}
\delta &=& \underset{\underset{\sigma(p) = p+q -1}{\sigma \in
    \Sk_{p,q-1}}}{\sum} \e (\sigma) \\ & & \phantom{@@@@@} B(G(
Y_{p+q}, Y_{\sigma(p+1)}, \dots, Y_{\sigma(p+q-1)}),F(Y_{\sigma(1)},
\dots, Y_{\sigma(p-1)}, Y_{p+q-1}))
\end{eqnarray*}

Now commute, up to a sign, $G( Y_{p+q}, Y_{\sigma(p+1)}, \dots,
Y_{\sigma(p+q-1)})$ and $Y_{p+q-1}$ to obtain:
\begin{eqnarray*}
\delta &=& - \underset{\underset{\sigma(p) = p+q -1}{\sigma \in
    \Sk_{p,q-1}}}{\sum} \e (\sigma) \\ & & \phantom{@@@}
B(F(Y_{\sigma(1)}, \dots, Y_{\sigma(p-1)},G( Y_{p+q}, Y_{\sigma(p+1)},
\dots, Y_{\sigma(p+q-1)})), Y_{p+q-1})\\ &=& - (-1)^{p-1} (-1)^{q-1}
\underset{\underset{\sigma(p) = p+q -1}{\sigma \in \Sk_{p,q-1}}}{\sum}
\e (\sigma) \\ & & \phantom{@@@} B(F(G( Y_{\sigma(p+1)}, \dots,
Y_{\sigma(p+q-1)}, Y_{p+q}), Y_{\sigma(1)}, \dots, Y_{\sigma(p-1)}),
Y_{p+q-1})
\end{eqnarray*}

Let $Z_i = Y_i$, $i = 1, \dots, p+q-2$ and $Z_{p+q-1} = Y_{p+q}$,
then:
\begin{eqnarray*}
F(G( Y_{\sigma(p+1)}, &\dots&, Y_{\sigma(p+q-1)}, Y_{p+q}),
Y_{\sigma(1)}, \dots, Y_{\sigma(p-1)}) \\ &=& F(G( Z_{\sigma(p+1)},
\dots, Z_{\sigma(p+q-1)}, Z_{p+q-1}), Z_{\sigma(1)}, \dots,
Z_{\sigma(p-1)}) \\ &=& F(G( Z_{\tau(1)}, \dots, Z_{\tau(q)}),
Z_{\tau(q+1)}, \dots, Z_{\tau(p+q-1)})
\end{eqnarray*}

where $\tau(1) = \sigma(p+1)$, $\dots$, $\tau(q-1) = \sigma(p+q-1)$,
$\tau(q) = p+q-1 = \sigma(p)$, $\tau(q+1) = \sigma(1)$, $\dots$,
$\tau(q+p-1) = \sigma(p-1)$. Comparing the inversions of $\tau$ with
the inversions of $\sigma$, it is easy to check that $\e (\tau) =
(-1)^{p-1} (-1)^{q-1} (-1)^{(p-1)(q-1)} \e(\sigma)$.

Finally
\begin{eqnarray*}
\delta &=& - (-1)^{(p-1)(q-1)} \underset{\underset{\tau(q) = p+q
    -1}{\tau \in \Sk_{q,p-1}}}{\sum} \e (\tau) \\ & & \phantom{@@@@}
B(F(G( Z_{\tau(1)}, \dots, Z_{\tau(q)}), Z_{\tau(q+1)}, \dots,
Z_{\tau(p+q-1)}), Y_{p+q-1})
\end{eqnarray*}

Then
\begin{eqnarray*}
\alpha + \delta &=& - (-1)^{(p-1)(q-1)} \underset{\tau \in
  \Sk_{q,p-1}}{\sum} \e (\tau) \\ & & \phantom{@@@@} B(F(G(
Z_{\tau(1)}, \dots, Z_{\tau(q)}), Z_{\tau(q+1)}, \dots,
Z_{\tau(p+q-1)}), Y_{p+q-1})
\end{eqnarray*}

Using similar arguments to compute $\beta + \gamma$, one obtains the
required identity.
\end{enumerate}

\end{proof}

\begin{rem}
When $\g$ is finite dimensional, there is a direct proof of
(\ref{5.1.4})(3) (avoiding computations) that we shall give in the
proof of Proposition \ref{6.4.1}, in Remark \ref{6.4.2}.
\end{rem}

\subsection{}

We assume now that $\g$ is a quadratic Lie algebra.

\begin{prop}\label{5.4.3}
$(\Cc_c(\g),d)$ is a subcomplex of the adjoint cohomology complex
  $(\Mcc_a(\g),d)$ of $\g$.
\end{prop}

\begin{proof}
It is enough to check that $d(\Cc_c(\g)) \subset \Cc_c(\g)$, but this
is obvious from Proposition \ref{5.1.4}(3) because $d = \ad(F_0)$ and
$F_0 \in \Cc_c(\g)$ since $\g$ is quadratic.
\end{proof}

\begin{defn}\label{5.4.4}
The cohomology of the complex $(\Cc_c(\g), d)$ is called the {\em
  cyclic cohomology} of $\g$, and denoted by $H^*_c(\g)$.
\end{defn}

\begin{rem}
Since $d = \ad (F_0)$, the Gerstenhaber bracket induces a \gla
structure on $H_c^*(\g)$.
\end{rem}

\begin{prop}\label{5.3.2bis}
The map $\Theta \colon \Cc_c(\g) \to \extgp$ is a homomorphism of
complexes. Moreover, $\Theta$ induces a map $\Theta^* \colon H_c^*(\g)
\to H^*_+(\g)$, which is an isomorphism when $\g$ is finite
dimensional.
\end{prop}

\begin{proof}
By an easy computation, one has $\Theta \circ d = \partial \circ
\Theta$ and the two first claims follow. For the third claim, use
Proposition \ref{5.1.4}.
\end{proof}

\begin{ex}\label{5.3.5bis}
  Assume that $\g$ is the Lie algebra associated to an associative
  algebra with a trace such that the bilinear form $\Tr(XY) := XY$,
  $\forall X$, $Y \in \g$ is non degenerate (e.g. $\g$ is the Lie
  algebra of finite rank operators on a given vector space, see
  Examples \ref{s.6} and \ref{s.7}). Consider the standard polynomials
  $\Ac_k$, for $k \geq 0$ if $\g$ had a unit, or for $k > 0$, if $\g$
  has no unit. Since $[\Ac_2, \Ac_{2k}]_a = 0$ by Proposition
  \ref{s.2.2}, each $\Ac_{2k}$ is a cocycle, then by Proposition
  \ref{s.5.7}, it is a cyclic cocycle, and one has $\Theta(\Ac_{2k}) =
  \frac{1}{2k+1} \ \Tr(\Ac_{2k+1})$.
\end{ex}

\subsection{} \label{5.4}

We assume now that $\g$ is a $n$-dimensional quadratic vector
space. Using the super Poisson bracket, we shall now go further into
the structure of $\Cc_c(\g)$. We need to renormalize the map $\Theta$,
defining $\Phi := - \frac12 \Theta$. We denote by $\mu$ the canonical
map from $\extg$ onto $\extgq$, and by $\Psi$ the map $\Psi = \mu
\circ \Phi$ from $\Cc_c(\g)$ into $\extgq$.

\begin{prop}\label{6.4.1} \hfill

\begin{enumerate}

\item If $C \in \Cc_c(\g)$, one has $\Dt(C) = \adpo(\Phi(C))$.

\item The restriction map $\Hht = \Dt|_{\Cc_c(\g)}$ is a \gla
  isomorphism from $\Cc_c(\g)[1]$ onto $\Hc(\g)$.

\item $\Psi$ is a \gla isomorphism from $\Cc_c(\g)[1]$ onto
  $\extgq[2]$.
\end{enumerate}

\end{prop}

\begin{proof}
  Fix an orthonormal basis $\{X_1, \dots, X_n \}$ of $\g$ and $\{
  \omega_1, \dots, \omega_n \}$ the dual basis.  Given $C \in
  \Cc_c(\g)$, $Y_1, \dots, Y_p \in \g$, 
\begin{eqnarray*}
  \adpo(\Phi(C))(\omega_k) (Y_1, &\dots&, Y_p) \\ &=& 2 (-1)^p \left(
  \sum_{r=1}^n \iota_{X_r} (\Phi(C)) \wedge \iota_{X_r} (\omega_r)
  \right) (Y_1, \dots, Y_p) \\ &=& (-1)^{p+1} B(C(X_k, Y_1, \dots,
  Y_{p-1}), Y_p) \\ &=& B(C(Y_1, \dots, Y_{p-1},X_k),Y_p) \\ &=& -
  B(C(Y_1, \dots, Y_{p-1},Y_p),X_k) \\ & =& - \omega_k (C(Y_1, \dots,
  Y_{p})) = - \Dt( C) (\omega_k) (Y_1, \dots, Y_p)
\end{eqnarray*}

by a formula given in \ref{2.3}, and this proves (1). From (1), we
deduce that $\Dt$ maps $\Cc_c(\g)$ into $\Hc(\g)$.

To prove (2), we remark that $\adpo \circ \ \Phi$ is onto by
Proposition \ref{5.1.4} (2), so $\Hht$ is onto, one to one and a \gla
homomorphism by Proposition \ref{2.1.4} and this proves (2).

\begin{rem}\label{6.4.2}
Let us give a direct proof of \ref{5.1.4} (2): given $C$, $C' \in
\Cc_c(\g)$, from the preceding results, we can assume that $C =
\Ft(\adpo(\Omega))$, $C'= \Ft(\adpo(\Omega'))$, with $\Omega$,
$\Omega' \in \extgp$. Then:
\[ [C,C']_a = [\Ft(\adpo(\Omega)), \Ft ( \adpo(\Omega'))]_a =
   \Ft \left( [\adpo(\Omega), \adpo(\Omega')] \right) = \Ft \left(
   \adpo(\{\Omega, \Omega'\}) \right) \qed \]
\end{rem}

To prove (3), we use the \gla isomorphism $\nu \colon \extgq \to
\Hc(\g)$ defined from $\adpo \colon \extg \to \Hc(\g)$, so one has
$\nu \left( \mu(\Omega) \right) = \adpo(\Omega)$, $\Omega \in \extg$,
and then $\nu \left( \Psi(C) \right) = \adpo (\Phi(C)) = \Hht(C)$,
$\forall \ C \in \Cc_c(\g)$, so $\Psi = \nu^{-1} \circ \ \Hht$.
\end{proof}

\begin{cor}
The \gla $\Cc_c(\g)$ is isomorphic to $\Hc(\g)$, and to
$\tilde{H}(n)$.
\end{cor}

Using $\Phi$, we can pull back the $\wedge$-product of $\Wedge_+ \g$
on $\Cc_c(\g)$ defining:

\begin{defn} \label{5.2.1} \hfill
\[C \wedge C' := \Phi^{-1} \left( \Phi(C)  \wedge \Phi(C') \right),
\forall \ C, C' \in \Cc_c(\g).\]
\end{defn}

Hence $\Cc_c(\g)$ becomes an associative algebra (without unit),
graded by $\Cc_c(\g)[-1]$. To describe the $\wedge$-product of
$\Cc_c(\g)$, we define a natural $\extg$-module structure on
$\Mcc_a(\g)$ by:
\[ \Omega \cdot  ( \alpha\otimes X) := (\Omega \wedge \alpha) \otimes
X, \forall \ \Omega, \alpha \in \extg, \ X \in \g\]

\begin{prop}\label{5.2.3}
If $C\in \Cc_c^k(\g)$, $C' \in \Cc_c^{k'}(\g)$, then $C \wedge C' \in
\Cc_c^{k+k'+1}(\g)$, and one has:
\[ C \wedge C' =  \Phi(C) \cdot C' + (-1)^{(k+1)(k'+1)} \Phi(C') \cdot
C.\] 
\end{prop}

\begin{proof}
Let $C'' = \Phi( C) \cdot C' + (-1)^{(k+1)(k'+1)} \Phi(C') \cdot
C$. Then
\begin{eqnarray*}
&&\Phi(C'') (Y_1, \dots, Y_{k+k'+2}) = \\&& \sum_{\sigma \in
    \Sk_{k+1,k'}} \e(\sigma) \Phi(C) (Y_{\sigma(1)}, \dots,
  Y_{\sigma(k+1)}) \Phi(C') ((Y_{\sigma(k+2)}, \dots,
  Y_{\sigma(k+k'+1)}, Y_{k+k'+2}) + \\ && \phantom{@@@@@@@}
  (-1)^{(k+1)(k'+1)} \\ && \sum_{\sigma \in \Sk_{k'+1,k}} \e(\sigma)
  \Phi(C') (Y_{\sigma(1)}, \dots, Y_{\sigma(k'+1)}) \Phi(C)
  ((Y_{\sigma(k'+2)}, \dots, Y_{\sigma(k+k'+1)}, Y_{k+k'+2}).
\end{eqnarray*}

In the first term of the right hand side, for each $\sigma$ define
$\tau$ by $\tau(i) = \sigma(i)$, $i \leq k + k' +1$, and $\tau(k+k'+2)
= k+k' +2$. In the second term, for each $\sigma$ define $\tau$ by
$\tau(1) = \sigma(k'+2)$, $\tau(k) = \sigma(k+k'+1)$, $\tau(k+1) =
k+k'+2$, $\tau(k+2) = \sigma(1)$, $\tau(k+3) = \sigma(2)$, $\dots$,
$\tau(k+k'+2) = \sigma(k' +1)$, then $\e(\tau) = (-1)^{(k+1)(k'+1)}
\e(\sigma)$, and one has:
\begin{eqnarray*}
&&\Phi(C'') (Y_1, \dots, Y_{k+k'+2}) = \\&&
  \underset{\underset{\tau(k+k'+2)=k+k'+2}{\tau \in
      \Sk_{k+1,k'+1}}}{\sum} \e(\tau) \Phi(C) (Y_{\tau(1)}, \dots,
  Y_{\tau(k+1)}) \Phi(C') ((Y_{\tau(k+2)}, \dots, Y_{\tau(k+k'+2)}) +
  \\&& \underset{\underset{\tau(k+1)=k+k'+2}{\tau \in
      \Sk_{k+1,k'+1}}}{\sum} \e(\tau) \Phi(C) (Y_{\tau(1)}, \dots,
  Y_{\tau(k+1)}) \Phi(C') ((Y_{\tau(k+2)}, \dots, Y_{\tau(k+k'+2)}) =
  \\ && \Phi(C) \wedge \Phi(C') (Y_1, \dots, Y_{k+k'+2}).
\end{eqnarray*}
\end{proof}

One has to be careful that $\ad(C)$ ($C \in \Cc_c(\g)$) is generally
not a derivation of the $\wedge$-product of $\Cc_c(\g)$, so the
following result is of interest:

\begin{prop} \label{5.4.6}
If $C\in \Cc_c^k(\g)$, $C' \in \Cc_c^{k'}(\g)$, $C'' \in \Cc_c(\g)$,
with $k \geq 1$, then:
\[ \ad(C)(C' \wedge C'') =  \ad(C)(C') \wedge C'' + (-1)^{(k+1)(k'+1)}
C'\wedge \ad(C)(C'').\]
\end{prop}

This means that when $C\in \Cc_c^k(\g)[1]$, then $\ad(C)$ is a
derivation of degree $k$ of the graded algebra $\Cc_c(\g)[-1]$ with
the $\wedge$-product.

\begin{proof}
One has 
\begin{eqnarray*}
\mu \left( \Phi( [C, C']_a) \right) &=& \Psi([C,C']) = [\Psi(C),
  \Psi(C')]_Q = [\mu \left( \Phi(C) \right), \mu \left( \Phi(C')
  \right)]_Q \\ &=& \mu \left( \{ \Phi(C), \Phi(C') \} \right)
\end{eqnarray*}
Since $\adpo(\Phi(C)) \left( \extg \right) \subset \extgp$, it follows
that $\Phi(\ad(C)(C')) = \adpo(\Phi(C)) \left( \Phi(C') \right)$, and
the result is proved using the fact that $\adpo(\Phi(C))$ is a
derivation of degree $k-1$ of $\extg$, and the definition of the
$\wedge$-product of $\Cc_c(\g)$.
\end{proof}

Using the Proposition above, and $d = \ad(F_0)$ with $F_0 \in
\Cc_c^2(\g)$, it results that the $\wedge$-product of $\Cc_c(\g)$
induces a $\wedge$-product on $H_c^*(\g)$ and $\phi^* = - \frac12
\theta^*$ is clearly an isomorphism of graded algebras from
$H^*_c(\g)$ onto $H^*_+(\g)$. From the definition of the \gla bracket
on $H^*_c(\g)$, denoting by $\mu^*$ the canonical map from $H^*(\g)$
onto $H_Q^*(\g)= H^*(\g) / \CC$, the map $\Psi^* = \mu^* \circ
\ \Phi^*$ is a \gla isomorphism from $H^*_c(\g)$ onto $H^*_Q(\g)$. We
summarize in:

\begin{prop}\label{5.4.6bis}
As a graded associative algebra, $H^*_c(\g)$ is isomorphic to
$H^*_+(\g)$ and as a \gla, $H_c^*(\g)$ is isomorphic to $H^*_Q(\g)$.
\end{prop}

\begin{ex}\label{5.4.7bis}
Let $\g = \gl(n)$. Then $H^*_+(\g) = \ext_+[a_1, a_3, \dots,
  a_{2n-1}]$, where $a_k = \Tr(\Ac_{k})$, $k \geq 0$
(e.g. \cite{Fuchs}). One has $\Theta(\Ac_{2k}) = \frac{1}{2k+1}
\ \Tr(\Ac_{2k+1})$ (Example \ref{5.3.5bis}), so by Proposition
\ref{5.4.6bis}, $H^*_c(\g) = \ext_+[\Ac_0, \Ac_2, \dots,
  \Ac_{2n-2}]$. The \gla bracket will be computed in Example
\ref{8.3.1}.
\end{ex}

\begin{rem}
When $\g$ is not finite dimensional, the map $\Theta^*$ of Proposition
\ref{5.3.2bis} is no longer an isomorphism, as shown with the
following example: let $V$ be an infinite dimensional vector space,
and $\g$ be the quadratic Lie algebra of finite rank operators of $V$,
as defined in Example \ref{s.7}. Recall that the invariant bilinear
form is $B(X,Y) = \Tr(XY)$, $X$, $Y \in \g$. Notice that $B(X,Y)$ is
well defined when $X \in \g$ and $Y \in \End(V)$. Moreover, the
formula $B([X,Y], Z) = - B(Y,[X,Z])$ is valid if at least one argument
is in $\g$. By Remark \ref{s.9}, $H^0_c(\g) = Z(\g) = \{0\}$ and
$H^1(\g) = \CC \ \Tr$, so:
\end{rem}

\begin{prop}\label{6.4.10}
  The map $\Theta^* \colon H^0_c(\g) \to H^1(\g)$ is not onto.
\end{prop}

Moreover,

\begin{prop}\label{6.4.11}
The map $\Theta^* \colon H^1_c(\g) \to H^2(\g)$ is not one to one.
\end{prop}

\begin{proof}
Fix $U \in \End(V)$ such that $U \notin \g \oplus \CC\ \Id_V$ and
consider the skew symmetric derivation $D$ of $ \g$ defined by $D =
\ad(U)|_\g$. The derivation $D$ is a cyclic cocycle but $D = \ad(Y)$
with $Y \in \g$ cannot be true because if $U' \in \End(V)$ commutes
with $\g$, then $U'$ must be a multiple of $\Id_V$. So $D$ is not a
cyclic coboundary. On the other hand, $\hat{D} (X,Y) = B(D(X),Y) =
\partial \omega(X,Y)$ where $\omega \in \Wedge^1 \g$ is defined by
$\omega(X) = - B(U,X)$, $X \in \g$. Hence $\hat{D}$ is a coboundary,
and if we denote by $\overline{D}$ the class of $D$ in $H_c^*(\g)$, we
get $\Theta^*(\overline{D}) =0$, and $\overline{D} \neq 0$.
\end{proof}

\section{The case of reductive and semisimple Lie
  algebras}\label{Section6}

\subsection{} \label{6.1}

Let $\g$ be a $n$-dimensional quadratic Lie algebra with bilinear form
$B$. We recall the natural $\g$-modules structures on $\extg$ and
$\Mcc_a(\g)$ defined by:
\[ \theta_X(\Omega)(Y_1, \dots, Y_p) = - \sum_i \Omega(Y_1, \dots, [X,
  Y_i], \dots, Y_p), \ \forall \ X, Y_1, \dots, Y_p \in \g, \Omega \in
\Wedge{}^p \g.\]
\[L_X(\Omega \otimes Y) = \theta_X(\Omega) \otimes Y + \Omega \otimes
   [X,Y], \ \forall X, Y \in \g, \Omega \in \extg.\]

Using the notation in \ref{5.4}, it is easy to check that
\[ \Phi \circ L_X = \theta_X \circ \Phi, \ \forall \ X \in
\g.\]

So one has:

\begin{prop}\label{6.1.4}
$\Cc_c(\g)$ is a $\g$-submodule of the $\g$-module $\Mcc_a(\g)$ and
  the isomorphism $\Phi$ (of \ref{5.4}) is a $\g$-module isomorphism
  from $\Cc_c(\g)$ onto $\Wedge_+ \g$.
\end{prop}

It is well known that any element of $\extgg$ is a cocycle, and if
$\g$ is reductive, that $H^\star(\g) = \extgg$ \cite{Koszul1}. Using
Propositions \ref{5.3.2bis}, \ref{5.4.6bis} and \ref{6.1.4}, we
deduce:

\begin{prop}\label{6.1.5}
Any invariant cyclic cochain is a cocycle. If $\g$ is reductive, any
cyclic cohomology class contains one, and only one invariant cyclic
cocycle (for instance, the only invariant cyclic coboundary is $0$).
\end{prop}

Hence, when $\g$ is reductive, we can identify $H_c^*(\g)$ and
$\Cc_c(\g)^\g$. This identification is valid for the corresponding
$\wedge$-products (actually isomorphic to the $\wedge$-product of
$(\Wedge_+ \g)^\g \simeq H_+^\star(\g)$) and for the corresponding
graded Lie bracket induced by the Gerstenhaber bracket (actually
isomorphic to $\Hc(\g)^\g$ and $\left( \extgq \right)^\g$).

\subsection{} \label{6.2}

In the remaining of this Section, we assume that $\g$ is a semisimple
Lie algebra with invariant bilinear form $B$ (not necessarily the
Killing form).

\begin{prop} \label{6.2.1} 
If $I$ and $I' \in \extgg$, then $\{I, I'\} = 0$.
\end{prop}

As a consequence of this Proposition and of Proposition
\ref{5.4.6bis}, one has:

\begin{cor} \label{6.2.2}
The Gerstenhaber bracket induces the null bracket on $H_c^*(\g) \simeq
\Cc_c(\g)^\g$.
\end{cor}

To prove Proposition \ref{6.2.1}, we need several lemmas: first, let
$\hk$ be a Lie algebra and $I \in \left( \Wedge^{p+1}
\hk\right)^\hk$. Define a map $\Omega \colon \hk \to \Wedge^{p} \hk$
by $\Omega(X) = \iota_X (I), \forall \ X \in \hk$. Then since $I$ is
invariant, one has:

\begin{lem}\label{6.2.3}
$\Omega$ is a morphism of $\hk$-modules from $(\hk, \ad)$ into
  $(\Wedge^p \hk, \theta)$.
\end{lem}

\begin{proof}
For all $X$, $Y$ and $Z \in \g$, we have: 
\[ \theta_X (\Omega(Y)) = \theta_X (\iota_Y (I)) = [\theta_X, \iota_Y] (I) +
\iota_Y(\theta_X(I)) = \iota_{[X,Y]} (I) = \Omega ([X,Y]). \]
\end{proof}

As a second argument for the proof of Proposition \ref{6.2.1}:

\begin{lem}\label{6.2.4}
Assuming that $\hk$ is a perfect Lie algebra (i.e. $\hk = [\hk,
  \hk]$), there exists a map $\alpha \colon \hk \to \Wedge^{p-1} \hk$
such that $\Omega = \partial \circ \alpha$ ($\partial$ is the
differential of the trivial cohomology of $\hk$). Moreover, if $\hk$
is semisimple, there exist an $\hk$-homomorphism $\alpha$ such that
$\Omega = \partial \circ \alpha$.
\end{lem}

\begin{proof}
If $X \in \hk$, we can find $Z_i$, $T_i \in \hk$ such that $X = \sum_i
[Z_i, T_i]$. Then, $\Omega(X) = \sum_i \theta_{Z_i} (\Omega(T_i))$ by
Lemma \ref{6.2.3}. But $\partial (\Omega(T_i)) = \partial
(\iota_{T_i}(I)) = \theta_{T_i} (I) - \iota_{T_i} (\partial(I)) =0$
since $I$ is an invariant. But $\theta_{Z_i}$ maps $Z^p(\hk)$ into
$B^p(\hk)$, so $\Omega(X) \in B^p (\hk)$. To construct $\alpha$, fix a
section $\sigma$ of the map $\partial \colon \Wedge^{p-1} \hk \to
B^p(\hk)$, i.e. $\sigma \colon B^p(\hk) \to \Wedge^{p-1} \hk$ such
that $\partial \circ \sigma = \Id_{B^p(\hk)}$ and then set $\alpha =
\sigma \circ \Omega$. When $\g$ is semisimple, one can fix a section
$\sigma$ which is a $\g$-homomorphism.
\end{proof}

\begin{proof}{\em (of Proposition \ref{6.2.1})} \hfill

Fix an orthonormal basis $\{X_1, \dots, X_n\}$ of $\g$ with respect to
$B$. Given $I$, $I' \in \extgg$, let $\Omega_r = \iota_{X_r} (I)$,
$\Omega'_r = \iota_{X_r}(I')$, $\alpha$, $\alpha'$ the
$\g$-homomorphisms given by Lemma \ref{6.2.4} and finally $\alpha_r=
\alpha(X_r)$, $\alpha_r' = \alpha'(X_r)$ so that $\Omega_r = \partial
\ \alpha_r$ and $\Omega_r' = \partial \ \alpha_r'$. With these
notations, in order to finish the proof, we need to show that $\sum_r
\Omega_r \wedge \Omega_r' = 0$. But:
\[ \sum_r \Omega_r \wedge \Omega_r' = \sum_r \partial \alpha_r \wedge
\partial \alpha_r' = \partial ( \sum_r \alpha_r \wedge \partial
\alpha_r') = 0\] 

since $\sum_r \alpha_r \wedge \partial \alpha_r' \in \extgg$.
\end{proof}

\begin{rem}\label{6.2.5}
Proposition \ref{6.2.1} can be directly deduced from a deep result of
Kostant \cite{Kostant2} about the structure of $\Cliff(\gd)^\g$ seen
as a deformation of $\extgg$: by the Hopf-Koszul-Samelson theorem,
$\extgg$ is an exterior algebra $\ext[a_1, \dots, a_r]$ with
$\rank(\g)$ $ = r$ and $a_1, \dots, a_r$ primitive (odd)
invariants. Kostant shows that $\Cliff(\gd)^\g$ is a Clifford algebra
constructed on $a_1, \dots, a_r$. Since the deformation from $\extg$
to $\Cliff(\gd)$ has leading term the Poisson bracket, it results that
$\{ a_i, a_j \} = 0$, $\forall \ i, j$, and then Proposition
\ref{6.2.1} follows.
\end{rem}

\begin{ex} \label{6.2.6}
Using the results in Section \ref{Section5}, and Corollary
\ref{6.2.2}, we will describe $H^*_c(\sk)$ and $H_c^*(\g)$ when $\sk=
\slk(n)$ and $\g = \gl(n)$ both equipped with the bilinear form
$B(X,Y) = \Tr(XY)$, $\forall \ X, Y$. Let $\mathbf{1}_\g$ be the
identity matrix.

One has $\Wedge \sk = \{ \Omega \in \Wedge \g \mid
\iota_{\mathbf{1}_\g}(\Omega) = 0 \}$ and $\Mcc_a(\sk) = \{ F \in
\Mcc_a(\g) \mid \iota_{\mathbf{1}_\g}(F) = 0 \text{ and } F(\g^p)
\subset \sk (F \in \Mcc_a^p(\g)) \}$. By Propositions \ref{s.4.2} and
\ref{s.3.6}, $\Ac_{2k} \in \Mcc_a(\sk)$. Moreover, let $a_k =
\Tr(\Ac_k)$ ($k \geq 0$), then by Proposition \ref{s.3.6}, $a_{2k+1}
\in \left( \extg \right)^\g$, $ \forall \ k \geq 0$, and by
Proposition \ref{s.4.2}, $a_{2k+1} \in \left( \Wedge \sk \right)^\sk$,
$ \forall \ k \geq 0$.

\begin{enumerate}

\item It is well known that $H^*(\g) = \extgg$ is the exterior algebra
  generated by the invariant cocycles $a_1$, $a_3$, $\dots$,
  $a_{2n-1}$, i.e. $\extgg = \ext[a_1,a_3, \dots, a_{2n-1}]$ and that
  $H^*(\sk) = \left( \Wedge \sk \right)^\sk$ is the exterior algebra
  generated by the invariant cocycles $a_3$, $a_5$, $\dots$,
  $a_{2n-1}$, i.e. $\left( \Wedge \sk \right)^\sk = \ext[a_3, a_5,
    \dots, a_{2n-1}]$ (see \cite{Kostant2, Kostant1, Fuchs}).

\item We need to compute the super Poisson bracket on $\extgg$. Note
  that \linebreak $\{ \Omega, \Omega' \} = 0$, $\forall \ \Omega,
  \Omega' \in \left( \Wedge \sk \right)^\sk$ by Proposition
  \ref{6.2.1}. Then, using $\sk^\perp = \CC \ \mathbf{1}_\g$, an
  adapted orthonormal basis, and the formula in Proposition
  \ref{4.2.1}, one finds that $\{ a_1, a_1 \} = 2n$. Then, since any
  element in $\extgg$ decomposes as $\Omega + \Omega' \wedge a_1$,
  with $\Omega$, $\Omega' \in \left( \Wedge \sk \right)^\sk$, we have
  only to compute the following brackets:
\begin{eqnarray*}
  && \{ \Omega, \Omega' \wedge a_1 \} = 0, \forall
  \ \Omega, \Omega' \in \ext^{w'}[a_3, \dots, a_{2n-1}], \\
  && \{ \Omega \wedge a_1, \Omega' \wedge a_1 \} = 2n (-1)^{w'}
  \ \Omega \wedge \Omega', \\ && \phantom{@@@@@@@@@@@@} \forall
  \ \Omega \in \ext[a_3, \dots, a_{2n-1}], \Omega' \in \ext^{w'}[a_3,
  \dots, a_{2n-1}].
\end{eqnarray*}

\item Use the isomorphism $\Phi^*$ of Proposition \ref{5.4.6bis} to
  find $H_c^*(\sk) = \Cc_c(\sk)^\sk$ and $H_c^*(\g) =
  \Cc_c(\g)^\g$. One has $[\Ac_2, \Ac_{2k}] = 0$ by Proposition
  \ref{s.2.2}, so $\Ac_{2k}$ is a cocycle, obviously
  $\g$-invariant. By Proposition \ref{s.5.7}, it is a cyclic cocycle,
  and $\Phi(\Ac_{2k}) = - \frac{1}{2(2k+1)} a_{2k+1}$. It results that
 \[ \qquad  H_c^*(\sk) = \ext_+[\Ac_2, \Ac_4, \dots, \Ac_{2n-2}] \text{ and }
 H_c^*(\g) = \ext_+[\Ac_0, \Ac_2, \dots, \Ac_{2n-2}].\]

\item Now we compute the Gerstenhaber bracket. For $H_c^*(\sk)$, by
  Corollary \ref{6.2.2}, the Gerstenhaber bracket vanishes. For
  $H_c^*(\g)$, we use the isomorphism $\Psi^*$ (see Proposition
  \ref{5.4.6bis}) combined with \ref{6.2.6} (3) and the commutation
  rules in $H^* (\g)$ computed in \ref{6.2.6} (2) from which the
  commutation rules in $H_Q^*(\g) = H^*(\g) / \CC$ are
  deduced. Finally the result is the following:
\begin{eqnarray*}
  && [F,F']_a = 0, \forall \ F,F' \in \ext_+[\Ac_2, \dots,
    \Ac_{2n-2}], \\ 
&& [\Ac_0,F]_a = 0, \forall \ F \in \ext_+[\Ac_0, \Ac_2, \dots,
    \Ac_{2n-2}], \\ 
&& [F,F' \wedge \Ac_0]_a = 0, \forall \ F,F' \in \ext_+[\Ac_2, \Ac_4,
    \dots, \Ac_{2n-2}], \\ 
&& [\Ac_0, F' \wedge \Ac_0]_a = \frac{n}{2} (-1)^{f'} \ F', \forall
  \  F' \in \ext_+^{f'}[\Ac_2, \Ac_4, \dots, \Ac_{2n-2}], \\ 
&& [F\wedge \Ac_0, F' \wedge \Ac_0]_a = \frac{n}{2} (-1)^{f'} \ F
  \wedge F', \\ 
&& \phantom{@@@} \forall \ F \in \ext_+[\Ac_2, \Ac_4, \dots,
    \Ac_{2n-2}], F' \in \ext_+^{f'}[\Ac_2, \Ac_4, \dots, \Ac_{2n-2}], 
\end{eqnarray*}
Remark that for the last result, one uses: $F' \in \ext_+^{f'}[\Ac_2,
  \Ac_4, \dots, \Ac_{2n-2}] \cap \Cc_c^{p'} (\g)$, then $p' = f' + 1$
mod 2 and $\Phi(F') \in \Wedge^{f'} \g$.

\end{enumerate}

\end{ex}

\section{Quadratic $2k$-Lie algebras and cyclic
  cochains}\label{Section7} 

\subsection{} \label{7.1}

Let $\g$ be a finite dimensional quadratic vector space with bilinear
form $B$. Given $D \in \Dc^{2k-1}$, $k \geq 1$ denote by $F = \Ft_D$
the associated (even) structure on $\g$ (see Sections \ref{Section1}
and \ref{Section2}), that we also denote by a bracket notation:
\[ [Y_1, \dots, Y_{2k}]  = F(Y_1, \dots, Y_{2k}), \forall \ Y_1, \dots,
Y_{2k} \in \g. \]

\begin{defn}\label{7.1.1}
The bilinear form $B$ is $F$-{\em invariant} (or $F$ is a {\em
  quadratic structure with bilinear form} $B$) if $B([Y_1, \dots,
  Y_{2k-1},Y],Z) = $ $- B(Y, [Y_1, \dots, Y_{2k-1},Z])$, $\forall
\ Y_1, \dots,$ $ Y_{2k-1}$, $Y$, $Z \in \g$.
\end{defn}

We introduce the linear maps $\ad_{Y_1, \dots, Y_{2k-1}} \colon \g \to
\g$ by:
\[ \ad_{Y_1, \dots, Y_{2k-1}}(Y)= [Y_1, \dots, Y_{2k-1},Y], \ \forall
\ Y_1, \dots, Y_{2k-1},Y \in \g.\]

It is obvious that

\begin{prop}\label{7.1.2}
  The bilinear form $B$ is $F$-invariant if and only if $\ad_{Y_1,
    \dots, Y_{2k-1}} \in {\mathfrak{o}}(B)$, $\forall \ Y_1, \dots,
  Y_{2k-1} \in \g$.
\end{prop}

The next Proposition results directly from Propositions \ref{5.1.4}
and \ref{6.4.1}.

\begin{prop}\label{7.1.3}
\hfill

\begin{enumerate}
\item $F$ is quadratic if and only if it is a cyclic cochain.

\item $F$ is quadratic if and only if there exists $I \in
  \Wedge^{2k+1} \g$ such that $D = - \frac{1}{2} \ad_{\tt{P}} (I)$ and
  in that case, one has $I(Y_1, \dots, Y_{2k+1}) = B ([Y_1, \dots,
  Y_{2k}], Y_{2k+1})$, $\forall \ Y_1, \dots,$ $Y_{2k+1} \in \g$.
\end{enumerate}

\end{prop}

\subsection{} \label{7.2}

Keeping the notations of Proposition \ref{7.1.3}, a quadratic $F$ will
define a $2k$-Lie algebra structure on $\g$ (namely a quadratic
$2k$-Lie algebra) if and only if:
\begin{equation} \label{7.1.4}
[F,F]_a =0 \quad \text{ or } \quad [D,D]=0 \quad \text{ or } \quad
\{I,I \} =0.
\end{equation}

Examples of quadratic $2k$-Lie algebras can be directly deduced from
Proposition \ref{6.2.1}: let us assume in the remaining of \ref{7.2},
that $\g$ is a semisimple Lie algebra with bilinear form $B$ (not
necessarily the Killing form). Then one has:

\begin{prop}\label{7.2.1}
  Any invariant even cyclic cochain in $\Mcc_a(\g)$ defines a
  quadratic $2k$-Lie algebra on $\g$.
\end{prop}

These examples were introduced for the first time in \cite{Per}, in
the case of primitive elements in $\extgg$ (we shall come back to the
construction in \cite{Per} later in this Section).

Let $F$ be an invariant even cyclic cochain, denote by:
\[ [Y_1, \dots, Y_{2k}]  = F(Y_1, \dots, Y_{2k}), \forall \ Y_1, \dots,
Y_{2k} \in \g \] 

the associated quadratic $2k$-bracket on $\g$. Let us introduce, as in
\ref{7.1}:
\[I([Y_1, \dots, Y_{2k+1}) = B([Y_1, \dots, Y_{2k}],Y_{2k+1}), \
\forall \ Y_1, \dots, Y_{2k+1} \in \g,\] 

and the associated derivation $D = -\frac{1}{2} \ad_{\tt{P}} (I)$ of
$\extg$. Since $[D,D] = 2 D^2 =0$, we can define the associated
cohomology on $\extg$ by
\begin{equation*}
H^\star(F) = Z(D)/B(D)
\end{equation*}
where $Z(D) = \ker(D)$ and $B(D) = \im(D)$. 

The following Lemma has to be compared with Formula \ref{2.2.2} of
\ref{2.1}:

\begin{prop}\label{7.2.3}
Let $\{ X_1, \dots, X_n \}$ be an orthonormal basis of $\g$ with
respect to $B$. Then there exist $\beta_1, \dots, \beta_n \in
\Wedge^{2k-1} \g$ such that:
\[ D = \frac{1}{2} \sum_r \beta_r \wedge \theta_{X_r}.\]
\end{prop}

\begin{proof}
Let $\{ \omega_1, \dots, \omega_n \}$ be the dual basis of $\{ X_1,
\dots, X_n \}$. One has $\theta_{X_r} (\omega_s)(Y) = B([X_r,X_s],Y)$
for all $Y \in \g$. So $\theta_{X_r} (\omega_s) = - \theta_{X_s}
(\omega_r)$ for all $r$, $s$. Define $\Omega(X) = \iota_X(I)$, $X \in
\g$. By Lemma \ref{6.2.4}, there exists a $\g$-homomorphism $\alpha
\colon \g \to \Wedge^{2k-1} \g$ such that $\Omega = \partial \circ
\alpha$. Define $\alpha_r = \alpha(X_r)$, then $\theta_{X_r}
(\alpha_s) = \alpha ([X_r,X_s])$, so one has $\theta_{X_r} (\alpha_s)
= - \theta_{X_s} (\alpha_r)$. Define $\Omega_r = \Omega(X_r) =
\partial \alpha_r$, the one has:
\[ D = - \frac{1}{2} \adpo (I) = - \sum_r \Omega_r \wedge \iota_{X_r}.\]
So $D(\omega_r) = - \partial \alpha_r$. Then using $\partial =
\frac{1}{2} \sum_s \omega_s \wedge \theta_{X_s}$ (\cite{Koszul1}), one
has:
\[ \partial \alpha_r = - \frac{1}{2} \sum_s \omega_s \wedge
\theta_{X_r} (\alpha_s) = - \frac{1}{2} \left( \theta_{X_r} (\sum_s
\omega_s \wedge \alpha_s) - \sum_s \theta_{X_r} (\omega_s) \wedge
\alpha_s \right).\] But $\sum_s \omega_s \wedge \alpha_s$ is
$\g$-invariant, so :
\[  \partial \alpha_r =  \frac{1}{2} \sum_s \theta_{X_r} (\omega_s) \wedge
\alpha_s = - \frac{1}{2} \sum_s \alpha_s \wedge \theta_{X_r}(\omega_s)
= \dfrac12 \sum_s \alpha_s \wedge \theta_{X_s} (\omega_r).\]

Therefore, since $D$ and $\sum_s \alpha_s \wedge \theta_{X_s}$ are
derivations of $\extg$, one has $D = - \frac12 \sum_s \alpha_s \wedge
\theta_{X_s}$, and setting $\beta_s = - \alpha_s$, the Proposition is
proved.
\end{proof}

From Proposition \ref{7.2.3}, we deduce:

\begin{prop}
One has $\extgg \subset Z(D)$.
\end{prop}

From the fact that $I \in \extgg$, $D$ is a $\g$-homomorphism of the
$\g$-module $\extg$, which is semisimple. By standard arguments
(\cite{Koszul1}), one deduces:

\begin{prop}
One has $\extgg \subset H^\star(F)$.
\end{prop}

When $F$ is the Lie algebra structure of $\g$, it is well known that
$H^\star(F) = \extgg$ (\cite{Koszul1}).

\subsection{} \label{7.3}

Let us now place the constructions in \cite{Per} in our context. We
assume that $\g$ is a semisimple Lie algebra of rank $r$ and fix a non
degenerate symmetric bilinear form $B$ (not necessarily Killing) on
$\g$. Let $\So (\g) = \sym(\gd)$. Using Chevalley's theorem, there
exist homogeneous invariants $Q_1, \dots, Q_r$ with $q_i = \deg(Q_i)$
such that $\So(\g)^\g = \CC[Q_1, \dots, Q_r]$. Let $t \colon
\So(\g)^\g \to \extgg$ be the Cartan-Chevalley transgression operator
(\cite{Cartan}, \cite{Chevalley}). By the Hopf-Koszul-Samelson theorem
(\cite{Cartan}, \cite{Chevalley}, \cite{Kostant2}), one has $\extgg =
\ext[t(Q_1), \dots, t(Q_r)]$ and $\deg(t(Q_i))= 2 q_i -1$. By
(\ref{7.1.4}) and Proposition \ref{6.2.1}, any odd element $I$ in
$\extgg$ defines a quadratic $2k$-Lie algebra structure on $\g$ (and
corresponding generalized Poisson bracket on $\gd$). As a particular
case, this works for $t(Q_i)$, $i=1, \dots, r$ which define a
$(2q_i-2)$-Lie algebra structure on $\g$ and a $\GPB$ on $\gd$, and
these are exactly the examples given in \cite{Per}, though in these
papers there are no citations, neither to Chevalley \cite{Chevalley},
nor to Cartan \cite{Cartan}. Let us insist that not only primitive
invariants (as sometimes claimed in \cite{Per}), but actually all odd
invariants do define $2k$-Lie algebra structures on $\g$ (Propositions
\ref{7.1.3} and \ref{7.2.1}).

\begin{ex}\label{8.3.1}
  Using the notation and the results of Example \ref{6.2.6}, let us
  consider the case of $\g = \gl(n)$, with bilinear form $B(X,Y) =
  \Tr(XY)$, $\forall X,Y \in \g$.  Consider $C = F + F' \wedge \Ac_0$
  with $F$, $F' \in \ext_+[\Ac_2, \Ac_4, \dots, \Ac_{2n-2}]$. In order
  to have $C$ an even element of $\Mcc_a(\g)$, we have to assume that
  $F \in \ext_+^{\text{odd}}[\Ac_2, \Ac_4, \dots, \Ac_{2n-2}]$ and $F'
  \in \ext_+^{\text{even}}[\Ac_2,\dots, \Ac_{2n-2}]$ (see the last
  remark in Example \ref{6.2.6} (4)). Moreover, we have to assume that
  $F$ and $F' \wedge \Ac_0$ have the same degree in $\Mcc_a(\g)$, say
  $2k$. Then, from commutation rules in \ref{6.2.6} (4), $C$ defines a
  $2k$-Lie algebra structure on $\g$ if and only if $F' \wedge F' =
  0$. This last condition is obviously satisfied if $F'$ is
  decomposable. For instance, if $n \geq 3$, $\alpha \Ac_8 + \beta
  \Ac_0 \wedge \Ac_2 \wedge \Ac_4$, $\alpha$, $\beta \in \CC$, defines
  a 8-Lie algebra structure on $\g$; if $n \geq 4$, $\alpha \Ac_{14} +
  \beta \Ac_0 \wedge \Ac_4 \wedge \Ac_8$, $\alpha$, $\beta \in \CC$,
  defines a 14-Lie algebra structure on $\g$.
\end{ex}

\end{document}